\documentclass[12pt]{amsart}
\usepackage{cancel}
\usepackage{amssymb}
\usepackage{hyperref}
\usepackage{enumitem}
\usepackage{amsmath}
\usepackage{latexsym}
\usepackage{extarrows}
\usepackage{enumitem}
\usepackage{txfonts}
\usepackage{mathtools}
\usepackage{bm}
\usepackage{tikz-cd}
\usepackage{xcolor}
\usepackage{quiver}
\usepackage{comment}
\usepackage{graphics}
\usepackage{float}
\usepackage{relsize}
\usepackage{bm}

\usepackage{cleveref}
\makeatletter
\@namedef{subjclassname@2020}{%
  \textup{2020} Mathematics Subject Classification}
\makeatother
\usepackage[T1]{fontenc}

\hypersetup{
    hidelinks,
    linktoc=page,
    }

\counterwithin{figure}{section}
\numberwithin{equation}{section}

\theoremstyle{plain}
\newtheorem{maintheorem}{Theorem}
\newtheorem{theorem}{Theorem}[section]
\newtheorem{proposition}[theorem]{Proposition}

\newtheorem{lemma}[theorem]{Lemma}
\newtheorem{corollary}[theorem]{Corollary}
\newtheorem{maincorollary}[maintheorem]{Corollary}

\theoremstyle{plain}
\newtheorem*{namedthm}{\namedthmname}
\newcounter{namedthm}
\makeatletter
	\newenvironment{named}[2]
	{\def\namedthmname{#1}
	\refstepcounter{namedthm}
	\namedthm[#2]\def\@currentlabel{#1}}
	{\endnamedthm}
\makeatother

\theoremstyle{definition}
\newtheorem{definition}[theorem]{Definition}

\newtheorem{remark}[theorem]{Remark}

\newtheorem{example}[theorem]{Example}

\newtheorem{question}[subsection]{Question}

\newcommand\E{\mathbb{E}}
\newcommand\Z{\mathbb{Z}}

\newcommand\Q{\mathbb{Q}}
\newcommand\T{\mathbb{T}}
\newcommand\C{\mathbb{C}}
\newcommand\I{\mathcal{I}}

\newcommand\N{\mathbb{N}}

\newcommand\M{\mathcal{M}}

\newcommand\XX{\mathrm{\bf X}}

\def\acts{\curvearrowright}

\newcommand\eps{\varepsilon}

\newcommand\Poly{\mathrm{Poly}}

\newcommand{\HK}{\mathrm{HK}}

\newcommand{\diag}{\mathrm{diag}}
\newcommand{\gen}{\textup{gen}}

\newcommand{\bracket}[1]{\left\llbracket #1 \right\rrbracket}
\newcommand{\Cor}[1]{Q_*^{\bracket{#1}}}
\newcommand{\supp}{\textup{supp}}
\newcommand{\diam}{\textup{diam}}

\begin{document}


\baselineskip=17pt


\title[Unique ergodicity and density finite sums in abelian groups]{Unique ergodicity of systems of finite order and the density finite sums theorem in arbitrary abelian groups}

\author[E. Ackelsberg]{Ethan Ackelsberg}
\address{Institute of Mathematics\\ \'{E}cole Polytechnique F\'{e}d\'{e}rale de Lausanne (EPFL) \\ Lausanne \\ 1015 Switzerland}
\email{ethan.ackelsberg@epfl.ch}

\author[A. Jamneshan]{Asgar Jamneshan}
\address{Faculty of Mathematics, TU Dresden, 01062 Dresden, Germany}
\email{asgar.jamneshan@tu-dresden.de}

\author[O. Shalom]{Or Shalom}
\address{Department of Mathematics\\ Bar Ilan University \\ 
Ramat Gan \\
5290002, Israel}
\email{Or.Shalom@math.biu.ac.il}

\date{\today}

\begin{abstract} 
Let $\Gamma$ be a countable abelian group and $k\geq 1$. We establish that every minimal topological $\Gamma$-system of order $k$ is uniquely ergodic. This result, which was previously only known for finitely generated groups for arbitrary $k$, answers a question of Glasner, Gutman, and Ye and yields the answer to multiple questions by Candela, Gonz\'alez-S\'anchez, and Szegedy. 
As applications of our dynamical results, we obtain the following combinatorial theorems. Let $A\subseteq \Gamma$ be a set of positive upper Banach density. Then there are infinite sets $B_1,\ldots, B_k\subseteq \Gamma$ such that $B_1+\dots+B_k\subseteq A$. If $\Gamma$ satisfies a necessary divisibility assumption, then there are also $t\in \Gamma$ and an infinite set $B\subseteq \Gamma$ such that $A-t$ contains the set of all sums of $m$ pairwise distinct elements of $B$ for all $m\in\{1,\ldots,k\}$. For arbitrary $k$, these density finite sums theorems were only known in the integers previously.
\end{abstract}



\maketitle


\section{Introduction}

Hindman's theorem \cite{hindman1974} asserts that any finite coloring of $\N$ contains a monochromatic IP-set, that is, a set of natural numbers which contains all finite sums of pairwise distinct elements of some infinite set. 

A recent breakthrough \cite{sumsets} by Kra, Moreira, Richter, and Robertson establishes a density version of Hindman's theorem: 
\begin{theorem}\label{thm-kmrr}
    If $A\subseteq \N$ has positive upper Banach density, then for every $k\geq 1$ there are an infinite set $B\subseteq \N$ and $t\in \N$ such that $A-t$ contains the set of all sums of $m$ pairwise distinct elements of $B$ for all $m\in\{1,\ldots,k\}$. 
\end{theorem}
The shift by $t$ is necessary in general as can be seen by taking $A$ to be the odd numbers. An example attributed to Straus (see, e.g., \cite{hindman1979}) exhibits sets of arbitrary large density no translate of which contains an IP-set, making the restriction on $m$ necessary. 
Sumsets of pairwise distinct elements are sometimes named restricted sumsets in the literature. 
As an immediate consequence of \Cref{thm-kmrr}, one obtains:
\begin{corollary}\label{cor-kmrr}
    If $A\subseteq \N$ has positive upper Banach density, then for every $k\geq 1$ there are an infinite sets $B_1,\ldots, B_k\subseteq \N$ such that the (unrestricted) sumset $B_1+\cdots + B_k \subseteq A$. 
\end{corollary}

\Cref{cor-kmrr} was originally proved in earlier work of Kra, Moreira, Richter, and Robertson \cite{BFRRinfinitesums}.

In this paper, we prove the following generalizations of \Cref{thm-kmrr} and \Cref{cor-kmrr}. Let $\Gamma$ be a countable abelian group. The upper Banach density of a subset $A\subseteq \Gamma$ is defined by the formula 
\begin{equation*}
                d^*(A) \coloneqq \sup_{\Phi} \limsup_{N \to \infty} \frac{|A \cap \Phi_N|}{|\Phi_N|},
\end{equation*}
where the supremum is taken over all F{\o}lner sequences $\Phi$ in $\Gamma$. 

\begin{maintheorem}\label{thm-main}
Let $\Gamma$ be a countably infinite abelian group, and let $k\geq 1$. Then the following are equivalent. 
\begin{itemize}
    \item[(i)] For every subset $A\subseteq \Gamma$ of positive upper Banach density, there exist an infinite set $B\subseteq \Gamma$ and $t\in\Gamma$ such that $A-t$ contains the set of all sums of $m$ pairwise distinct elements of $B$ for all $m\in\{1,\ldots,k\}$. 
    \item[(ii)] $[\Gamma\colon m\Gamma]<\infty$ for all $m\in\{1,\ldots,k\}$. 
\end{itemize}
\end{maintheorem}
The implication $(i)\Rightarrow(ii)$ was established by the first author in \cite{sumset-counterexamples}. The main aim of this paper is to establish the converse implication. 

We obtain immediately an extension of \Cref{cor-kmrr} to countably infinite abelian groups under the divisibility assumption (ii) in \Cref{thm-main}. Our second combinatorial result proves a density unrestricted sumset result without this divisibility assumption.

\begin{maintheorem}\label{thm-main1}
Let $\Gamma$ be an arbitrary countably infinite abelian group, let $k\geq 1$, and let $A\subseteq \Gamma$ be a set of positive upper Banach density. Then there exist infinite sets $B_1,\ldots, B_k\subseteq \Gamma$ such that the (unrestricted) sumset $B_1+\cdots + B_k \subseteq A$. 
\end{maintheorem}
As a consequence, we obtain a multiplicative version of \Cref{cor-kmrr}:
\begin{maincorollary} \label{cor:main}
    Let $A \subseteq \N$ be a set of positive multiplicative upper Banach density.
    Then for every $k \in \N$, there exist infinite sets $B_1, \ldots, B_k \subseteq \N$ such that
    \begin{equation*}
        B_1 \cdot \ldots \cdot B_k = \left\{ \prod_{j=1}^k b_j : b_j \in B_j \right\} \subseteq A.
    \end{equation*}
\end{maincorollary}

\subsection*{Overview of previous work}
Moreira, Richter, and Robertson \cite{mrr2019} resolved the $k=2$ case of \Cref{cor-kmrr}; more generally, they established the corresponding product-set theorem for subsets of infinite amenable groups having positive upper density along a two-sided F\o lner sequence. Their proof combines ultrafilter techniques with a Jacobs--de Leeuw--Glicksberg-type decomposition. Host \cite{host2019sumsets} subsequently gave a shorter ergodic proof of their result, using the Kronecker factor and a suitable self-joining of an associated ergodic system. 
Kra, Moreira, Richter, and Robertson \cite{BFRRinfinitesums,BFRRb+b+t} developed a new dynamical approach to these questions using Host--Kra structure theory \cite{host2005nonconventional}, encoding infinite sumsets through recurrent configurations in self-joinings, which they called Erd\H{o}s cubes and Erd\H{o}s progressions.  This approach yielded \Cref{cor-kmrr} \cite{BFRRinfinitesums} and led to a proof of \Cref{thm-kmrr}, first for $k=2$ \cite{BFRRb+b+t} and then for arbitrary $k$ \cite{sumsets}.

Their methods provided a starting point for subsequent extensions to other ambient groups and more general configurations. Charamaras and Mountakis \cite{cm} established the $k=2$ case of \Cref{thm-main}, and obtained analogous results for ordered products in several classes of amenable groups. The first two authors \cite{AJ} proved the $k=3$ case of \Cref{thm-main} by relying on an ergodic structure theorem for systems of order~2 \cite{jst} due to the last two authors and Tao. A key new ingredient in the proof in \cite{AJ} was an equidistribution theorem for translational systems formed on degree $2$ filtered Polish groups.  

Further refinements, variants, and obstructions such as density criteria for unrestricted sumsets, asymmetric and polynomial configurations, limitations on the growth of the summands, and alternative proofs are studied in \cite{sumset-counterexamples,kr-unrestricted,kousek-asymmetric,ckmr-unrestricted,kousek-B+B,ackelsberg-rational,hh-growth,kmrr-short-bc,radic-uniform,hr-uniform,hkr-finite-sums,tz-sumset-primes}; see also the survey \cite{kmrr-survey-paper} for a broader account and further questions.

While putting the final touches on the current paper, a different, AI-assisted proof of \Cref{thm-main1} was independently obtained in \cite{li-su}.

\subsection*{Discussion of main results}
The first step in the proof of both \Cref{thm-main} and \Cref{thm-main1} is to translate the combinatorial problem to a dynamical one. Here, we closely follow the strategy developed in \cite{sumsets} (which underlies also its precedents \cite{BFRRinfinitesums,BFRRb+b+t}). Most of the necessary modifications to the more general setting of countable abelian groups are routine, though, at some point, we need to pass to a suitable topological model of an ergodic system of order $k$. Fortunately, for arbitrary abelian group actions  such a model was already constructed \cite{candela-szegedy-cubic} by Candela and Szegedy.

The remaining key dynamical problem is the content of the main new contribution of this work:  

\begin{maintheorem} \label{thm:minimal finite order systems are u.e.}
Let $\Gamma$ be an arbitrary countable abelian group, and let $k\geq 1$. Then every minimal topological $\Gamma$-system of order $k$ is uniquely ergodic. 
\end{maintheorem}

Our unique ergodicity result answers several questions raised in topological dynamics and nilspace theory which we address after discussing previously known special cases of it and our proof of it for arbitrary $k,\Gamma$.  

For arbitrary $k$, \Cref{thm:minimal finite order systems are u.e.} was previously only known when $\Gamma$ is a finitely generated group. 
Existing proofs in this case establish that minimal $\Gamma$-systems of order $k$ are $k$-step pro-nilsystems first. Then, since unique ergodicity is preserved under taking inverse limits, \Cref{thm:minimal finite order systems are u.e.} follows from Leibman's results \cite[Theorem 2.19 and Corollary 2.20]{leibman2005pointwise} on unique ergodicity of minimal $\Gamma$-nilsystems (the strategy goes back to Parry's work on nilflows \cite{parry1969ergodic}).
The aforementioned topological structure theorem was established for $\Z$-systems in the case $k=2$ by Host and Maass \cite{hm2007} and for arbitrary $k$ by Host, Kra, and Maass \cite{hkm} using Host--Kra theory. The topological structure theorem for groups admitting a dense compactly generated subgroup (such as finitely generated discrete groups) was established by Gutman, Manners, and Varj\'u \cite{gmv,gmv2,gmv3} using nilspace theory. 
The first two authors \cite{AJ} established the $k=2$ case of \Cref{thm:minimal finite order systems are u.e.} for arbitrary abelian $\Gamma$. To this end, \cite{AJ} proves a Ratner-type equidistribution theorem for translational $\Gamma$-systems formed on degree $2$ filtered Polish groups and deduces unique ergodicity for  minimal topological $\Gamma$-system of order $2$ from the ergodic structure theorem \cite{jst} of the last two authors and Tao classifying ergodic measure-preserving $\Gamma$-systems of order $2$ as inverse limits of such translational systems.

Unfortunately, a structure theorem analogous to the Host--Kra structure theorem for general $\Gamma$-systems of order $k$ is  currently not known and constitutes a difficult open problem in the subject. Accordingly, we take a substantially different approach to proving \Cref{thm:minimal finite order systems are u.e.}, which we now describe. 

Our starting point is a recent topological structure theorem for minimal $\Z^\omega$-systems of finite order \cite{c_g-s_sz} by Candela, Gonz\'alez--S\'anchez, and Szegedy that states that minimal $\Z^\omega$-systems of finite order are dynamical factors of inverse limits of polynomial orbit $\Z^\omega$-systems (see \Cref{def:polynomial-orbit-system}). We observe (\Cref{thm:polynomial-orbit-system}) that every polynomial orbit $\Z^\omega$-system is a minimal component of a translational $\Z^\omega$-system (\Cref{def:translational system}). Thus, it suffices to prove unique ergodicity for orbit closures in translational $\Gamma$-systems.

The ideas to prove unique ergodicity, originating in \cite{AJ}, differ from Leibman's approach due to the fact that nilpotent Polish groups comprise a much larger class than nilpotent Lie groups.
In particular, Leibman's approach (inspired by Parry's earlier arguments) proceeds by induction on dimension of the nilpotent Lie group. The nilpotent groups with which we work may not be finite-dimensional, so the strategy is to perform a complexity-reduction process on the degree of the system rather than its dimension.
Degree~1 translational systems are rotations on compact abelian groups. Unique ergodicity of abelian rotational systems is classical (see, e.g., \cite[Proposition 2.3]{AJ}). 
The induction step is broken into two smaller steps.
First, we show, using some ideas from Host--Kra structure theory, that the unique ergodicity statement for degree $k$ translational systems implies a Wiener--Wintner theorem for ergodic averages in general measure-preserving systems  weighted by continuous functions evaluated along orbits in degree $k$ translational systems. 
Then we adapt an argument of Lesigne \cite{lesigne} to show that this Wiener--Wintner theorem can be combined with the nilpotent structure of translational systems to produce unique ergodicity of orbits in translational systems of degree $k+1$. 

Thus, along the way, we also obtain a generalization of the nilsequence Wiener--Wintner theorem for $\Z^d$-actions (see, \cite[Theorem 2.22]{hkuniformity} by Host--Kra  for $\Z$-actions, and \cite[Theorem 2.4 and specification (i) after Theorem 2.4]{fhweighted} by Frantzikinakis--Host for $\Z^d$-actions) to the setting of arbitrary countable abelian group actions. 
Say that $\psi \colon \Gamma \to \C$ is a \emph{$\Gamma$-sequence of order $k$} if there exists a topological dynamical $\Gamma$-system $(Y,S)$ of order $k$, a continuous function $F \colon Y \to \C$, and a point $y \in Y$ such that $\psi(\gamma) = F(T^{\gamma}y)$.

\begin{maintheorem} \label{thm:WW}
    Let $\Gamma$ be a countable abelian group.
    Let $\Phi$ be a tempered F{\o}lner sequence in $\Gamma$.
    Suppose $(X, \mu, T)$ is a measure-preserving $\Gamma$-system.
    Then for every $f \in L^1(\mu)$, there exists a subset $X'_f \subseteq X$ with $\mu(X'_f) = 1$ such that if $x \in X'_f$ and $\psi \colon \Gamma \to \C$ is a $\Gamma$-sequence of order $k$, then
    \begin{equation*}
        \frac{1}{|\Phi_N|} \sum_{\gamma \in \Phi_N} \psi(\gamma) f(T^{\gamma}x)
    \end{equation*}
    converges.
\end{maintheorem}

 \Cref{thm:minimal finite order systems are u.e.} answers several questions from the literature. It answers affirmatively \cite[Question 10.3]{ggy} for countable discrete abelian groups, provides an affirmative answer to \cite[Question 5.2]{CGSSextensions} using \cite[Theorem 2.9]{CGSSextensions}, and since distality and transitivity give minimality, it also answers affirmatively \cite[Question 5.5]{CGSSextensions} for countable discrete abelian groups. 

\bigskip

 We close this introduction with some questions. 

\begin{question}\label{question}
    Let $(G/\Lambda, \phi)$ be a $\Z^{\omega}$-translational system (\Cref{def:translational system}).
    Is it true that all minimal components are of the form $Hx = \{hx : h \in H\}$ for some closed subgroup $H \leq G$ and $x \in G/\Lambda$?
\end{question}

Based on \cite[Theorem 1.5]{c_g-s_sz} and \Cref{thm:polynomial-orbit-system}, a positive answer to \Cref{question} would yield that any finite order minimal $\Z^\omega$-system is a factor of an inverse limit of translational $\Z^\omega$-systems. In particular, this would imply that ergodic Host--Kra factors of order $k$ for $\Z^{\omega}$-systems are factors of inverse limits of degree $k$ translational $\Z^\omega$-systems. 
The $k=2$ case of \Cref{question} was established in \cite{AJ}.

\begin{question}\label{question1}
    Let $\Gamma = (\Gamma,\cdot)$ be an infinite amenable group and $k\geq 1$.
    Is it true that for every set $A \subseteq \Gamma$ having positive density with respect to a two-sided F{\o}lner sequence, there exist infinite sets $B_1,\ldots,B_k\subseteq \Gamma$ such that $B_1\cdot B_2\cdot \ldots \cdot B_k\subseteq A$? 
\end{question}

As already mentioned, the $k=2$ case of \Cref{question1} was resolved affirmatively in \cite[Theorem 1.3]{mrr2019}.
The $k=3$ case was asked explicitly in \cite[Question 5.19]{kmrr-survey-paper}.
It is also an interesting question to know what happens if $A$ is only assumed to have positive density with respect to a one-sided F{\o}lner sequence (see \cite[Question 5.20]{kmrr-survey-paper}).

\subsection{Organization of the paper}
In \Cref{sec:preliminaries}, we define topological and measure-preserving systems of order $k$ and collect some of their basic properties. \Cref{sec:polynomial orbit systems} refines a result of Candela, Gonz\'alez-S\'anchez, and Szegedy, who established that any minimal $\Z^\omega$-system of finite order is a factor of an inverse limit of polynomial orbit systems. We show that these polynomial orbit systems can be embedded as orbit closures in translational systems (see \Cref{thm:polynomial-orbit-system}). \Cref{sec:unique ergodicity} proves that any such orbit closure is uniquely ergodic (see \Cref{thm:equidistribution}), which leads to the proof of \Cref{thm:minimal finite order systems are u.e.}. \Cref{sec:unique ergodicity} is also devoted to generalizing the Wiener--Wintner theorem to general abelian groups and a more general form of nilsequences (see \Cref{thm:WW-translational}); this implies \Cref{thm:WW}. The remaining sections apply the dynamical results to derive the combinatorial theorems. \Cref{sec:dynamical reformulation} formulates the combinatorial results of this paper as dynamical problems. \Cref{sec:proof of density FS} is then devoted to the proof of \Cref{thm-main}, \Cref{sec:proof sum of k sets} to \Cref{thm-main1}, and \Cref{sec:proof of corollary} to \Cref{cor:main}.

\subsection*{Acknowledgements}
E.A. was supported by the Swiss National Science Foundation under grant TMSGI2-211214.
A.J. was supported by the German Research Foundation under Germany's
Excellence Strategy -- EXC-2047 -- 390685813 and its Heisenberg Programme --
547294463. 
O.S. is supported by Alon Fellowship and Polak foundation. We would like to thank Pablo Candela for explaining the results in his paper with Szegedy \cite{candela-szegedy-cubic} pertaining to the topological model employed in this paper. We would like to thank Yonatan Gutman who answered some of our questions about his paper with Lian \cite{gutman-lian} which inspired some of our results.

\subsection*{AI Statement}
While writing the introduction, LLMs were used to locate references. 
No other LLM usage was involved.

\section{Preliminaries} \label{sec:preliminaries}

We assume familiarity with the basic definitions in topological dynamics and ergodic theory, including topological dynamical systems, measure-preserving systems, topological and measurable factors, and notions such as transitivity, minimality, distality, and ergodicity.

Some objects and notations we use in the paper are less standardized, so we record them in this preliminary section for reference. Throughout this section, assume that $\Gamma$ is an arbitrary countable discrete abelian group. 

\subsection{Generic points}

Let $(X, T)$ a topological dynamical $\Gamma$-system and  $\mu$ a $T$-invariant Borel probability measure on $X$.
Let $\Phi = (\Phi_N)_{N \in \N}$ be a F{\o}lner sequence in $\Gamma$. 
A point $x \in X$ is \emph{generic} for $\mu$ along $\Phi$, written $x \in \gen(\mu,\Phi)$, if
\begin{equation*}
    \lim_{N \to \infty} \frac{1}{|\Phi_N|} \sum_{\gamma \in \Phi_N} \delta_{T^{\gamma}x} = \mu
\end{equation*}
in the weak* topology.
That is, for every continuous function $f \in C(X)$,
\begin{equation*}
    \lim_{N \to \infty} \frac{1}{|\Phi_N|} \sum_{\gamma \in \Phi_N} f(T^{\gamma}x) = \int_X f~d\mu.
\end{equation*}

\subsection{Notations and conventions for discrete cubes}

For $k \geq 0$, we write $\bracket{k}$ for the $k$-dimensional discrete hypercube $\bracket{k} = \{0,1\}^k$.
We write elements of $\bracket{k}$ as strings $\bm{\omega} = \omega_1\ldots\omega_k$ and equip $\bracket{k}$ with the co-lexicographic order.
For example, the elements of $\bracket{2}$ written in order are $00, 10, 01, 11$.
Given an element $\bm{\omega} \in \bracket{k}$, the number of 1s appearing in $\omega$ is denoted by $|\bm{\omega}| = |\{i \in [k] : \omega_i = 1\}|$.

\subsection{Regionally proximal relations and topological systems of order $k$}

\begin{definition}[Order $k$ in the topological settings]
Let $(X,T)$ be a minimal topological $\Gamma$-dynamical system and $k\geq 1$.
    \begin{itemize}
        \item[(i)] The set of $k$-dimensional \emph{dynamical cubes} of $X$ is the set
        $$Q^{\bracket{k}}(X)\coloneqq \overline{\{(T^{\bm{\omega} \cdot \bm{\gamma}}x)_{\bm{\omega}\in\bracket{k}} :\bm{\gamma} \in \Gamma^k, x \in X \}}.$$
        \item[(ii)] We say that $x,y\in X$ are $(k-1)$-regionally proximal, and write $x\sim_{k-1}y$, if $(x,y,y,\dots,y)\in Q^{\bracket{k}}$. 
        \item[(iii)] A system $X$ is \emph{topologically of order $k$} if the $k$-regionally proximal relation is trivial, meaning $x\sim_k y\Longrightarrow x=y$.
        \item[(iv)] The set of \emph{$k$-dimensional dynamical corners} of $X$ is 
        $$\Cor{k}\coloneqq \{\bm{x}\in X^{\bracket{k}\backslash \{\bm{0}\}} : \exists\, y\in X \text{ such that } (y,\bm{x})\in Q^{\bracket{k}}\}.$$
    \end{itemize}
\end{definition}

A reference for the following facts in a general setting are \cite[Theorem 3.8 and Theorem 7.14]{ggy}. 

\begin{lemma}[Properties of the regional proximal relation] \label{lem: RP properties}
Let $(X,T)$ be a minimal topological $\Gamma$-dynamical system and $k\geq 1$.
\begin{itemize}
    \item[(i)] $\sim_k$ is a closed equivalence relation.
    \item[(ii)] $\Cor{k} \subseteq X^{\bracket{k}\backslash\{0^k\}}$ is compact.
    \item[(iii)] If $X$ is topologically of order $k$, then there exists a continuous map $\tau\colon\Cor{k+1} \rightarrow X$ such that $(\tau(\bm{x}),\bm{x})\in Q^{\bracket{k+1}}$ for all $\bm{x}\in \Cor{k+1}.$
\end{itemize}
     Moreover, a system is topologically of order $k$ if and only if there is unique $(k+1)$-cube completion, i.e. for every $\bm{x} \in \Cor{k+1}$, there exists a unique element $y \in X$ such that $(y,\bm{x}) \in Q^{\bracket{k+1}}$.
\end{lemma}

\begin{definition}[Parallelepiped group] \label{def:parallelepiped group}
    Let $(X, T)$ be a topological dynamical $\Gamma$-system and $k \geq 1$.
    \begin{itemize}
        \item The \emph{diagonal action} is the action $T_{\Delta} \colon \Gamma \acts Q^{\bracket{k}}$ defined by $(T_{\Delta}^{\gamma} \bm{x})_{\bm{\omega}} = T^{\gamma}x_{\bm{\omega}}$ for every $\bm{\omega} \in \bracket{k}$.
        \item For $i \in [k]$, the \emph{face action} is the action $T_i \colon \Gamma \acts Q^{\bracket{k}}$ defined by
        \begin{equation*}
            (T_i^{\gamma} \bm{x})_{\bm{\omega}} = \begin{cases}
                T^{\gamma} x_{\bm{\omega}}, & \text{if}~\omega_i = 1; \\
                x_{\bm{\omega}}, & \text{if}~\omega_i=0.
            \end{cases}
        \end{equation*}
        \item The parallelepiped group $\mathcal{G}^{\bracket{k}}$ is the group generated by the diagonal and face transformations. 
    \end{itemize}
\end{definition}

\begin{lemma}[{\cite[Lemma 4.1]{hkm},\cite[Proposition 4.6 and Lemma 7.11]{ggy}}] \label{lem: minimal distal cube system}
    Let $(X, T)$ be a topological $\Gamma$-system.
    If $(X,T)$ is minimal and distal, then $(Q^{\bracket{k}},\mathcal{G}^{\bracket{k}})$ is minimal and distal.
\end{lemma}

\subsection{Host--Kra factors and measure-preserving systems of order $k$}

Let $(X,\mu,T)$ be a measure-preserving $\Gamma$-system. We assume familiarity with Host--Kra theory. In particular, the Host--Kra seminorms $\|\cdot\|_{U^k}$, measures $\mu^{\bracket{k}}$, and factors $Z_k(X)$ are all defined in \cite{host2005nonconventional} (see also \cite{hk-book}). While \cite{host2005nonconventional,hk-book} develops Host--Kra theory only for $\Z$-systems, it was previously observed on several occasions (see, e.g., \cite[Appendix A]{btz}, \cite[\S 2]{jst}) that the basic theory directly extends to arbitrary countable abelian group actions.  

The following inequality is well known.

\begin{lemma}[Host--Kra--Gowers--Cauchy--Schwarz inequality {\cite[Lemma 3.9]{host2005nonconventional}}] \label{lem: hkgcs inequality}
    Let $(X,\mu,T)$ be an ergodic measure-preserving $\Gamma$-system and $k \in \N$.
    If $(f_{\bm{\omega}})_{\bm{\omega} \in \bracket{k}}$ is a family of bounded functions on $X$, then
    \begin{equation*}
        \left| \int_{X^{\bracket{k}}} \bigotimes_{\bm{\omega} \in \bracket{k}} f_{\bm{\omega}}~d\mu^{\bracket{k}} \right| \leq \prod_{\bm{\omega} \in \bracket{k}} \| f_{\bm{\omega}} \|_{U^k}
    \end{equation*}
\end{lemma}

\begin{proposition} \label{prop: topological order k implies measurable order k}
    Suppose $(X, T)$ is a minimal dynamical $\Gamma$-system of order $k$.
    If $\mu$ is an ergodic $T$-invariant Borel probability measure on $X$, then $(X, \mu, T)$ is measurably of order $k$.
\end{proposition}
\begin{proof}
    We begin with a couple of preliminary observations.
    
    Let $\rho \colon X^{\bracket{k+1}} \to X^{\bracket{k+1} \setminus \{\bm{0}\}}$ be the projection map $$\rho((x_{\bm{\omega}})_{\bm{\omega} \in \bracket{k+1}}) = (x_{\bm{\omega}})_{\bm{\omega} \in \bracket{k+1} \setminus \{\bm{0}\}},$$ and let $\mu_*^{\bracket{k+1}} = \rho_* \mu^{\bracket{k+1}}$.
    Since $\mu^{\bracket{k+1}}(Q^{\bracket{k+1}}) = 1$ (see \cite[Lemma 5.6]{hkm} or \cite[Lemma 5.6]{gutman-lian}, the same proof works for infinitely generated groups), we have $\mu_*^{\bracket{k+1}}(\Cor{k+1}) = 1$.

    By \Cref{lem: RP properties}, there exists a continuous function $\tau\colon\Cor{k+1} \rightarrow X$ such that $(\tau(\bm{x}),\bm{x})$ is a $(k+1)$-dimensional cube for all $\bm{x} \in \Cor{k+1}$.
    If $\bm{x} \in Q^{\bracket{k+1}}$, then the two points $\bm{x} = (x_{\bm{0}}, \rho(\bm{x}))$ and $(\tau(\rho(\bm{x})), \rho(\bm{x}))$ are both $(k+1)$-cubes, so $\tau(\rho(\bm{x})) = x_{\bm{0}}$ by the unique completion property of systems of order $k$.
    Therefore, $\tau_* \mu_*^{\bracket{k+1}} = (\tau \circ \rho)_* \mu^{\bracket{k+1}} = \mu$ by \cite[Lemma A.15]{btz}. 
    
    We now follow the strategy outlined in \cite[p. 303]{gmv}.

    Let $f \colon X \to \C$ be a bounded measurable function that is not almost everywhere equal to zero, say with $\|f\|_{L^2(\mu)} = \delta > 0$.
    We will prove that $\|f\|_{U^{k+1}}>0$.
   
   Define $F \colon X^{\bracket{k+1} \setminus \{\bm{0}\}} \to \C$ by
   \begin{equation*}
       F(\bm{x}) = \begin{cases}
           f(\tau(\bm{x})), & \text{if}~\bm{x} \in \Cor{k+1}; \\
           0, & \text{otherwise}.
       \end{cases}
   \end{equation*}
  Then we have
  \begin{equation} \begin{split} \label{eq: L^2 as correlation}
    \int_{X^{\bracket{k+1}}} f \otimes \overline{F}~d\mu^{\bracket{k+1}}
    & = \int_{Q^{\bracket{k+1}}} f \otimes \overline{F}~d\mu^{\bracket{k+1}} \\
    & = \int_{\Cor{k+1}} f(\tau(\bm{x})) \cdot \overline{F(\bm{x})}~d\mu_*^{\bracket{k+1}}(\bm{x}) \\
     & = \int_{\Cor{k+1}} (f \circ \tau) \cdot (\overline{f \circ \tau})~d\mu_*^{\bracket{k+1}} \\
     & = \int_X |f|^2~d\mu = \delta^2.
  \end{split} \end{equation}

  Now, since $\mu_*^{\bracket{k+1}}$ is a Borel measure on the product space $X^{\bracket{k+1}\backslash \{\bm{0}\}}$, there exist $\ell \in \N$ and a family of continuous functions $H^{(i)}_{\bm{\omega}} \in C(X)$, $i \in [\ell]$, $\bm{\omega} \in \bracket{k+1} \setminus \{\bm{0}\}$, such that
  \begin{equation} \label{eq: tensor representation}
      \left\| F - \sum_{i=1}^{\ell} \bigotimes_{\bm{\omega} \in \bracket{k+1} \setminus \{\bm{0}\}} H^{(i)}_{\bm{\omega}} \right\|_{L^2\left( \mu_*^{\bracket{k+1}} \right)} < \delta.
  \end{equation}
   Combining \eqref{eq: L^2 as correlation} and \eqref{eq: tensor representation}, there exists $i \in [\ell]$ such that
   $$\left| \int_{X^{\bracket{k+1}}} f(x_{\bm{0}})\cdot \prod_{\bm{\omega} \in \bracket{k+1} \setminus \{\bm{0}\}} \overline{H^{(i)}_{\bm{\omega}}(x_{\bm{\omega}})}~d\mu^{\bracket{k+1}} \right| >0.$$
   On the other hand, by the Host--Kra--Gowers--Cauchy--Schwarz inequality (\Cref{lem: hkgcs inequality}),
   $$\left| \int_{X^{\bracket{k+1}}} f(x_{\bm{0}})\cdot \prod_{\bm{\omega} \in \bracket{k+1} \setminus \{\bm{0}\}} \overline{H^{(i)}_{\bm{\omega}}(x_{\bm{\omega}})}~d\mu^{\bracket{k+1}} \right| \leq \|f\|_{U^{k+1}} \cdot \prod_{\bm{\omega} \in \bracket{k+1} \setminus \{\bm{0}\}} \|H^{(i)}_{\bm{\omega}}\|_{U^{k+1}},$$
   so $\|f\|_{U^{k+1}} > 0$, as required.
\end{proof}

\section{Polynomial orbit systems are subsystems of translational systems} \label{sec:polynomial orbit systems}

Our strategy to prove \Cref{thm:minimal finite order systems are u.e.} is to find a way to identify nilpotent structures in systems of order $k$. We fix our notation regarding these nilpotent structures in the following definition. 

\begin{definition}[Filtered groups and polynomial maps]\label{def:polynomial-orbit-system}
Let \(G\) be a topological\footnote{All topological groups considered in this paper are Hausdorff.} group equipped with a sequence $G_\bullet=(G_i)_{i=0}^{\infty}$ of closed subgroups.
\begin{itemize}
    \item
    We call $G_\bullet=(G_i)_{i=0}^{\infty}$ a \emph{prefiltration} on \(G\) if
    \[
    G_0 \supseteq G_1 \supseteq G_2 \supseteq \cdots
    \]
    and
    \[
    [G_i,G_j]\subseteq G_{i+j}
    \qquad\text{for all }i,j\geq 0.
    \]
    We say that \(G_\bullet\) has \emph{degree at most \(k\)}, or simply \emph{degree \(k\)}, if $G_{k+1}=\{1\}$.
    We refer to the pair \((G,G_\bullet)\) as a \emph{prefiltered group}. 
    \item For any $j\geq 0$, we denote by
    \[
    G_{\bullet+j}:=(G_{i+j})_{i=0}^{\infty}
    \]
    the \emph{shifted prefiltration}.
    \item     If additionally
    \[
    G_0=G_1=G,
    \]
    we say that \(G_\bullet\) is a \emph{filtration}, and that
    \((G,G_\bullet)\) is a \emph{filtered group}. 
    \item If $G$ admits a finite degree filtration, then it is \emph{nilpotent}.
    \item Given any closed subgroup \(\Lambda\leq G\), the \emph{induced prefiltration} on \(\Lambda\) is
    \[
    \Lambda_\bullet=(\Lambda_i)_{i=0}^{\infty},
    \qquad
    \Lambda_i:=G_i\cap\Lambda.
    \]
    If \(G_\bullet\) is a filtration, then \(\Lambda_\bullet\) is also a filtration.\footnote{Indeed,
    $[\Lambda_i,\Lambda_j]\subseteq \Lambda\cap [G_i,G_j]\subseteq \Lambda\cap G_{i+j}=\Lambda_{i+j}$.}

    We say that \(\Lambda\) is \emph{adapted to \(G_\bullet\)}, if $G_i\Lambda$ is closed in $G$ for every $i\geq0$.
    In this case we write
    \[
    (G/\Lambda)_\bullet:=(G_i/\Lambda_i)_{i=0}^{\infty}.
    \]

   \item A central class of examples of filtered groups are \emph{filtered nilmanifolds} $(G/\Lambda, G_\bullet)$, where \((G,G_\bullet)\) is a filtered Lie\footnote{We explicitly include disconnected Lie groups.} group of finite degree and \(\Lambda\leq G\) is a lattice adapted to $G_\bullet$.     
    \item 
    If \(G=(G,+)\) is an additive group, its \emph{abelian filtration} is
    \[
    G_0=G_1:=G,
    \qquad
    G_i:=\{0\}
    \quad\text{for }i>1.
    \]
    Throughout, we equip \(\Gamma=\Z^\omega\) with its abelian filtration.

    \item
    Let \(G=(G,G_\bullet)\) be a prefiltered group and let \(k\geq 0\).
    The \emph{Host--Kra group}
    \[
    \HK^k(G_\bullet)
    \leq G^{\bracket{k}}
    \]
    is the subgroup generated by the elements
    \[
    \left(
        g_{\omega_0}^{\,\mathbf{1}_{\bm{\omega} \geq \bm{\omega}_0}}
    \right)_{\bm{\omega}\in\bracket{k}},
    \]
    where $\bm{\omega}_0\in\bracket{k}$, $g_{\bm{\omega}_0}\in G_{|\bm{\omega}_0|}$, and \(\bm{\omega}\geq\bm{\omega}_0\) means that $\omega_j\geq(\omega_0)_j$ for every $j=1,\ldots,k$.

    \item
    Let \(G=(G,G_\bullet)\) be a prefiltered group and let \(\Lambda\leq G\) be a closed subgroup adapted to $G_\bullet$. If $\pi:G\to G/\Lambda$ denotes the quotient map, define
    \[
    \HK^k((G/\Lambda)_\bullet)
    :=
    \pi^{\bracket{k}}\bigl(\HK^k(G_\bullet)\bigr)
    \subseteq
    (G/\Lambda)^{\bracket{k}},
    \]
    where
    \[
    \pi^{\bracket{k}}
    \left((g_{\bm{\omega}})_{\bm{\omega} \in \bracket{k}}\right)
    :=
    \left(\pi(g_{\bm{\omega}})\right)_{\bm{\omega}\in\bracket{k}}.
    \]
    There is a canonical identification, as a homogeneous-space quotient,
    \[
    \HK^k((G/\Lambda)_\bullet)
    =
    \HK^k(G_\bullet)/\HK^k(\Lambda_\bullet),
    \]
    where $\Lambda_i:=\Lambda\cap G_i$.

    \item Let \((G,G_\bullet)\) be a prefiltered group and let \(\Lambda\leq G\) be a closed subgroup adapted to $G_\bullet$. A map $p:\Z^\omega\to G/\Lambda$ is called a \emph{polynomial map} if for every \(k\geq0\),
    \[
    p^{\bracket{k}}
    \bigl(
    \HK^k((\Z^\omega)_\bullet)
    \bigr)
    \subseteq
    \HK^k((G/\Lambda)_\bullet),
    \]
    where
    \[
    p^{\bracket{k}}
    \left(
    (x_{\bm{\omega}})_{\bm{\omega}\in\bracket{k}}
    \right)
    :=
    \left(
    p(x_{\bm{\omega}})
    \right)_{\bm{\omega}\in\bracket{k}}.
    \]
    We denote the set of all such polynomial maps by
    \[
    \Poly(\Z^\omega,(G/\Lambda)_\bullet).
    \]
    If \(\Lambda=\{1\}\), $p:\Z^\omega\to G$ is a polynomial map if and only if
    \begin{equation}\label{eq:derivative}
    \partial_{h_1}\cdots\partial_{h_k}p(x)\in G_k
    \end{equation}
    for every \(k\geq0\) and every $x,h_1,\ldots,h_k\in\Z^\omega$, where the difference operator \(\partial_h\) is defined by
    \[
    \partial_hp(x)
    :=
    p(x+h)p(x)^{-1}.
    \]
    \item We define the shift action $\sigma\colon \Z^\omega\times \Poly(\Z^\omega,G_\bullet) \to \Poly(\Z^\omega,G_\bullet)$ by 
\begin{equation}\label{eq:shift action}
    \sigma^\gamma(p)(\eta)\coloneqq p(\eta+\gamma),
    \qquad \eta,\gamma\in\Z^\omega.   
\end{equation}
\end{itemize}
\end{definition}

We collect some useful basic properties of polynomial maps. 

\begin{lemma}[Polynomial calculus]\label{lem:poly-calculus}
For all $i,j\geq0$ and $\gamma\in\Z^\omega$,
\[
[\Poly(\Z^\omega,G_{\bullet+i}),\Poly(\Z^\omega,G_{\bullet+j})]
\subseteq
\Poly(\Z^\omega,G_{\bullet+i+j}),
\]
\[
\sigma^\gamma\bigl(\Poly(\Z^\omega,G_{\bullet+i})\bigr)
=
\Poly(\Z^\omega,G_{\bullet+i}),
\]
and, for every $p\in\Poly(\Z^\omega,G_{\bullet+i})$,
\[
\sigma^\gamma(p)p^{-1}
\in
\Poly(\Z^\omega,G_{\bullet+i+1}).
\]
Moreover, each $\Poly(\Z^\omega,G_{\bullet+i})$ is closed in $G^{\Z^\omega}$.
\end{lemma}

\begin{proof}
    The commutator property follows from \cite[Proposition B.6]{gtz}, see the example after \cite[Corollary B.5]{gtz} for the translation property, and the third property follows from the derivative definition \eqref{eq:derivative} of a polynomial map. Since a prefiltration consists of closed subgroups and \eqref{eq:derivative} depends continuously on finitely many coordinates, each $\Poly(\Z^\omega,G_{\bullet+i})$ is closed in $G^{\Z^\omega}$. 
\end{proof}

In order to study the nilpotent structures in dynamical systems, we equip filtered groups with a dynamical structure as follows. 

\begin{definition}[Translational $\Gamma$-systems]
\label{def:translational system}
Let $G$ be a Polish group with a filtration $G_\bullet$ of degree $k$, and let $\Lambda\leq G$ be a closed subgroup adapted to $G_\bullet$.
Let $\Gamma$ be a countable abelian group and $\phi\colon\Gamma\to G$ a group homomorphism.
The topological dynamical $\Gamma$-system $(G/\Lambda,T_\phi)$, where
\[
T_\phi^\gamma x=\phi(\gamma)x,
\qquad x\in G/\Lambda,\quad \gamma\in\Gamma,
\]
is called a \emph{degree $k$ translational $\Gamma$-system}.
We will often abuse notation and write $(G/\Lambda,\phi)$ instead of $(G/\Lambda,T_\phi)$. 
\end{definition}

The starting point of our analysis is a recent structural result \cite[Theorem 1.5]{c_g-s_sz} by Candela, Gonz\'alez-S\'anchez and Szegedy which identifies minimal finite order topological $\Z^\omega$-systems with factors of inverse limits of so-called polynomial orbit systems which we now pause to define. 

\begin{definition}\label{def:polynomial orbit system}
Let $k\geq 1$. A topological $\Z^\omega$-system $(Y,S)$ is called a \emph{degree-$k$ polynomial orbit $\Z^\omega$-system} if there are a degree-$k$ filtered nilmanifold $(G/\Lambda,G_\bullet)$ and a polynomial map $p\in\Poly(\Z^\omega,(G/\Lambda)_\bullet)$ such that 
\[
Y=
\overline{\{\sigma^\gamma(p):\gamma\in\Z^\omega\}} \qquad \text{and} \qquad S=\sigma,
\]
where $\sigma$ is the shift action \eqref{eq:shift action}.   
\end{definition}

Polynomial orbit $\Z$-systems were studied in \cite[\S 14.2]{hk-book} by Host and Kra, who proved that they are $\Z$-nilsystems. 
In general, polynomial orbit $\Z^\omega$-systems fail to be $\Z^\omega$-nilsystems as the following example illustrates. 

\begin{example}
Consider $\Poly(\Z^\omega,\mathcal{D}^2(\T))$ where $\mathcal{D}^2(\T)$ is the degree $2$ filtration given by $G_0=G_1=G_2=\T$ and $G_3=0$. Every element in $\Poly(\Z^\omega,\mathcal{D}^2(\T))$ is a polynomial of degree at most $2$. These take the form $\gamma\mapsto b(\gamma,\gamma)+\phi(\gamma)+c$ where $b:\Z^\omega\times \Z^\omega\rightarrow \T$ is bilinear, $\phi:\Z^\omega\rightarrow \T$ is a character, and $c\in\T$. Using binomial coefficients and reindexing the quadratic coefficients, we have, as a group,
$$\Poly(\Z^\omega,\mathcal{D}^2(\T))\cong \T^{\N\times \N}\times \T^{\N}\times \T.$$
Now, let $b(\gamma,\gamma') = \sum_{i=1}^\infty \gamma_i\cdot\gamma_i'\cdot\alpha_i$ where all $\alpha_1,\alpha_2,\dots\in \T$ are irrational. Let $Q(\gamma) = b(\gamma,\gamma) = \sum_{i=1}^\infty \gamma_i^2\alpha_i$ and consider the polynomial orbit defined by $\mathcal{O}(Q)=\overline{\{Q\circ \sigma^\gamma : \gamma\in \Z^\omega\}}$. Then $(Q\circ \sigma^\gamma)(\gamma') = b(\gamma+\gamma',\gamma+\gamma') = Q(\gamma') + 2b(\gamma',\gamma) + b(\gamma,\gamma)$. Hence,
$$Q\circ \sigma^\gamma = Q+\phi_\gamma + c_\gamma$$
where $\phi_\gamma(\gamma') = 2\sum_{i=1}^\infty \gamma_i'\cdot\gamma_i\cdot\alpha_i$ is a character and $c_\gamma = \sum_{i=1}^\infty \gamma_i^2\alpha_i\in \T$. Hence, $\mathcal{O}(Q)\subseteq Q+\Poly(\Z^\omega,\mathcal{D}^1(\T))\cong \T^\N\times\T$, where $\mathcal{D}^1(\T)$ is the abelian filtration on $\T$. We claim that $\mathcal{O}(Q)$ is not finite dimensional. Let $\pi:\Poly(\Z^\omega,\mathcal{D}^2(\T))\rightarrow \Poly(\Z^\omega,\mathcal{D}^2(\T))/\T$ denote the quotient map by the constants. We first show that $\pi(\mathcal{O}(Q)) = \pi(Q+\Poly(\Z^\omega,\mathcal{D}^1(\T)))$. Since $\mathcal{O}(Q)$ is compact, the image $\pi(\mathcal{O}(Q))$ is closed and so it suffices to prove it is dense. Let $\phi\in \Poly(\Z^\omega,\mathcal{D}^1(\T))$ and choose a representative modulo constants with $\phi(0)=0$. Then $\phi$ is a character of $\Z^\omega$ and therefore takes the form $\phi(\gamma) = \sum_{i=1}^\infty \gamma_i\beta_i$ for some $\beta=(\beta_1,\beta_2,\dots)\in \T^\N$. Since all the coordinates of $\alpha$ are irrational, we can find a sequence $\gamma_n\in \Z^\omega$ such that $(2(\gamma_n)_i\alpha_i)_{i\geq 1}\rightarrow \beta$ in $\T^\N$. Therefore, modulo constants, $Q\circ \sigma^{\gamma_n}\rightarrow Q+\phi$, establishing the equality.

Finally, fix $t\in\T$. By Weyl's equidistribution theorem, we can choose $m_j\in\Z$ such that $m_j^2\alpha_j\to t$. Since every $\gamma'\in\Z^\omega$ has finite support, $Q(\gamma'+m_je_j)=Q(\gamma')+m_j^2\alpha_j$ for all sufficiently large $j$. Thus $Q+t\in\mathcal{O}(Q)$, and taking orbit closures gives $\mathcal{O}(Q)+t\subseteq\mathcal{O}(Q)$ for every $t\in\T$. The equality modulo constants now gives $\mathcal{O}(Q)=Q+\Poly(\Z^\omega,\mathcal{D}^1(\T))\cong\T^\N\times\T$. Thus $\mathcal{O}(Q)$ is infinite dimensional and the polynomial orbit system generated by $Q$ is not a $\Z^\omega$-nilsystem.
\end{example}

Adapting the arguments in \cite[\S 14.2]{hk-book}, our goal in this section is to show that polynomial $\Z^\omega$-orbit systems equivariantly embed into translational $\Z^\omega$-systems of the same degree (see \Cref{thm:polynomial-orbit-system}). We first realize the space of polynomial maps as a compact homogeneous space, and then pass to a closed subgroup whose filtration has degree at most $k$.

\begin{proposition}\label{prop:hat G}
Let $k\geq 1$ and let $(G/\Lambda,G_\bullet)$ be a degree-$k$ filtered nilmanifold.
\begin{itemize}
    \item[(i)] The canonical map induces a homeomorphism of compact spaces
    \[
    \Poly(\Z^\omega,G_\bullet)/\Poly(\Z^\omega,\Lambda_\bullet)
    \longrightarrow
    \Poly(\Z^\omega,(G/\Lambda)_\bullet).
    \]

    \item[(ii)] Using the shift action \eqref{eq:shift action}, define the semi-direct products 
    \[
    \hat G\coloneqq \Poly(\Z^\omega,G_\bullet)\rtimes_\sigma \Z^\omega,
    \qquad
    \hat \Lambda\coloneqq \Poly(\Z^\omega,\Lambda_\bullet)\rtimes_\sigma \Z^\omega.
    \]
    Then $\hat G$ is a Polish group and $\hat \Lambda$ is a closed cocompact subgroup of $\hat G$.
    Moreover, there is a homeomorphism
    \[
    \xi\colon \hat G/\hat \Lambda\to \Poly(\Z^\omega,(G/\Lambda)_\bullet),
    \qquad
    \xi((p,\gamma)\hat \Lambda)\coloneqq p\Lambda,
    \]
    satisfying for all $\gamma\in\Z^\omega$ and $x\in \hat G/\hat \Lambda$,
    \[
    \xi((1,\gamma)x)=\sigma^\gamma(\xi(x)).
    \]
\end{itemize}
\end{proposition}

\begin{proof}
By the Lazard--Leibman theorem (see, e.g., \cite[Proposition B.6]{gtz}), $\Poly(\Z^\omega,G_\bullet)$ is a group. Equip $\Poly(\Z^\omega,G_\bullet)$ with the subspace topology inherited from $G^{\Z^\omega}$. By \Cref{lem:poly-calculus}, this makes $\Poly(\Z^\omega,G_\bullet)$ a Polish group, and $\Poly(\Z^\omega,\Lambda_\bullet)$ is a closed subgroup.

Similarly, we equip $\Poly(\Z^\omega,(G/\Lambda)_\bullet)$ with the topology inherited from $(G/\Lambda)^{\Z^\omega}$. By Tychonoff's theorem and the construction of Host--Kra spaces, $\Poly(\Z^\omega,(G/\Lambda)_\bullet)$ is a compact Polish space.

By \cite[Lemma 5.5]{c_g-s_sz}, the map
\[
\Poly(\Z^\omega,G_\bullet)\to \Poly(\Z^\omega,(G/\Lambda)_\bullet),
\qquad p\mapsto p\Lambda,
\]
is surjective and, by construction, continuous. Moreover, $p\Lambda=q\Lambda$ if and only if
\[
q^{-1}p\in \Poly(\Z^\omega,\Lambda_\bullet).
\]
Thus the induced map
\[
\Poly(\Z^\omega,G_\bullet)/\Poly(\Z^\omega,\Lambda_\bullet)
\longrightarrow
\Poly(\Z^\omega,(G/\Lambda)_\bullet)
\]
is a continuous bijection. The action of $\Poly(\Z^\omega,G_\bullet)$ on the target by pointwise left multiplication is continuous and transitive, so by Effros's open mapping theorem (see, e.g., \cite[Theorem 22.1]{kechris1995classical}) this bijection is a homeomorphism. This proves (i).

We prove (ii). The group $\Z^\omega$, equipped with the discrete topology, is Polish, and the shift action $\sigma$ is continuous. Hence $\hat G$ is a Polish group and $\hat\Lambda$ is a closed subgroup. Define
\[
\Xi\colon \hat G/\hat \Lambda\to
\Poly(\Z^\omega,G_\bullet)/\Poly(\Z^\omega,\Lambda_\bullet),
\qquad
\Xi((p,\gamma)\hat \Lambda)
=
p\Poly(\Z^\omega,\Lambda_\bullet).
\]
This map is well-defined and bijective, since every $(p,\gamma)\hat \Lambda$ has a representative $(p,0)$ and
\[
(p,0)\hat \Lambda=(q,0)\hat \Lambda
\quad\Longleftrightarrow\quad
q^{-1}p\in \Poly(\Z^\omega,\Lambda_\bullet).
\]
The map $\Xi$ is continuous by the quotient topology, and its inverse
\[
p\Poly(\Z^\omega,\Lambda_\bullet)
\longmapsto
(p,0)\hat\Lambda
\]
is continuous for the same reason. Hence $\Xi$ is a homeomorphism.

Let $\xi$ be the composition of $\Xi$ with the homeomorphism from (i). In particular, $\hat\Lambda$ is cocompact in $\hat G$. Finally,
\[
\begin{aligned}
\xi((1,\gamma)(p,0)\hat \Lambda)
&=\xi((\sigma^\gamma(p),\gamma)\hat\Lambda)\\
&=\sigma^\gamma(p)\Lambda
=\sigma^\gamma(p\Lambda),
\end{aligned}
\]
which proves the equivariance.
\end{proof}

\begin{corollary}\label{cor:shifted-quotients}
For every $i\geq0$, the quotient
\[
\Poly(\Z^\omega,G_{\bullet+i})/
\Poly(\Z^\omega,\Lambda_{\bullet+i})
\]
is compact.
\end{corollary}

\begin{proof}
The proof of \Cref{prop:hat G}(i) applies verbatim to the shifted prefiltrations $G_{\bullet+i}$ and $\Lambda_{\bullet+i}$.
\end{proof}

The sequence 
\[
\Poly(\Z^\omega,G_{\bullet+i})\rtimes_\sigma\Z^\omega,
\qquad 2\leq i\leq k,
\]
generally fails to be a degree-$k$ filtration. Indeed, since $$(1,t)\in \Poly(\Z^\omega,G_{\bullet+j})\rtimes_\sigma\Z^\omega$$ for every $2\leq j\leq k$, if $p\in\Poly(\Z^\omega,G_{\bullet+i})$, then
\[
[(p,0),(1,t)]
=
\bigl(p\sigma^t(p^{-1}),0\bigr),
\]
which in general only lies in $\Poly(\Z^\omega,G_{\bullet+i+1})\times\{0\}$, whereas the filtration condition would require membership in $\Poly(\Z^\omega,G_{\bullet+i+j})\rtimes_\sigma\Z^\omega$. To overcome this issue, we adapt a construction in \cite[\S 14.2]{hk-book} due to Host and Kra: 

\begin{proposition}\label{prop:bar G}
Let $(G/\Lambda,G_\bullet)$ be a degree-$k$ filtered nilmanifold and put
\[
\bar G\coloneqq \Poly(\Z^\omega,G_{\bullet+1})\rtimes_\sigma \Z^\omega.
\]
Then $\bar G$ is a closed normal subgroup of $\hat G$. The sequence
\[
\bar G_0=\bar G_1\coloneqq \bar G,
\qquad
\bar G_i\coloneqq \Poly(\Z^\omega,G_{\bullet+i})\times\{0\},
\quad i\geq2,
\]
forms a filtration of degree at most $k$ on $\bar G$.
\end{proposition}

\begin{proof}
The subgroup $\Poly(\Z^\omega,G_{\bullet+1})$ is closed in $\Poly(\Z^\omega,G_\bullet)$ by \Cref{lem:poly-calculus}, hence $\bar G$ is closed in $\hat G$.

We prove normality. Let $(p,\gamma)\in \hat G$ and $(q,\eta)\in \bar G$. Then
\[
\begin{aligned}
(p,\gamma)(q,\eta)(p,\gamma)^{-1}
&=
\bigl(p\sigma^\gamma(q)\sigma^\eta(p^{-1}),\eta\bigr)\\
&=
\bigl(\bigl(p\sigma^\gamma(q)p^{-1}\bigr)
       \bigl(p\sigma^\eta(p^{-1})\bigr),\eta\bigr).
\end{aligned}
\]
By \Cref{lem:poly-calculus},
\[
p\sigma^\gamma(q)p^{-1}
=
[p,\sigma^\gamma(q)]\,\sigma^\gamma(q)
\in \Poly(\Z^\omega,G_{\bullet+1}),
\]
and
\[
p\sigma^\eta(p^{-1})
=
\bigl(\partial_\eta p\bigr)^{-1}
\in \Poly(\Z^\omega,G_{\bullet+1}).
\]
Thus
\[
(p,\gamma)(q,\eta)(p,\gamma)^{-1}\in\bar G.
\]

For $i\geq2$, the subgroup $\bar G_i$ is closed and normal in $\hat G$. Indeed, for $(p,\gamma)\in\hat G$ and $(q,0)\in\bar G_i$,
\[
(p,\gamma)(q,0)(p,\gamma)^{-1}
=
\bigl(p\sigma^\gamma(q)p^{-1},0\bigr)
\in\bar G_i
\]
by \Cref{lem:poly-calculus}.

We now verify the filtration property. The sequence is decreasing by definition. If $i,j\geq2$, then for $(p,0)\in\bar G_i$ and $(q,0)\in\bar G_j$,
\[
[(p,0),(q,0)]
=
([p,q],0)
\in\bar G_{i+j}.
\]
Next, let $i\geq2$. Since $\bar G$ is generated by
\[
\Poly(\Z^\omega,G_{\bullet+1})\times\{0\}
\]
and the elements $(1,\gamma)$, $\gamma\in\Z^\omega$, it is enough to note that, for $p\in\Poly(\Z^\omega,G_{\bullet+1})$ and $q\in\Poly(\Z^\omega,G_{\bullet+i})$,
\[
[(p,0),(q,0)]
=
([p,q],0)
\in\bar G_{i+1},
\]
and
\[
[(1,\gamma),(q,0)]
=
\bigl(\sigma^\gamma(q)q^{-1},0\bigr)
\in\bar G_{i+1}.
\]
Finally, if $p,q\in\Poly(\Z^\omega,G_{\bullet+1})$, then
\[
[(p,0),(q,0)]\in\bar G_2,
\qquad
[(1,\gamma),(q,0)]\in\bar G_2,
\]
and
\[
[(1,\gamma),(1,\eta)]
=
(1,0)\in\bar G_2.
\]
Hence $[\bar G,\bar G]\subseteq\bar G_2$. Since $\bar G_0=\bar G_1=\bar G$, these cases give
\[
[\bar G_i,\bar G_j]\subseteq\bar G_{i+j}
\qquad\text{for all }i,j\geq0.
\]
Since $G_\bullet$ has degree at most $k$,
\[
\bar G_{k+1}
=
\Poly(\Z^\omega,G_{\bullet+k+1})\times\{0\}
=
\{(1,0)\}.
\]
Thus $\bar G_\bullet$ has degree at most $k$.
\end{proof}

\begin{lemma}\label{lem:bar Lambda}
Put
\[
\bar \Lambda
\coloneqq
\Poly(\Z^\omega,\Lambda_{\bullet+1})\rtimes_\sigma\Z^\omega.
\]
Then
\[
\bar\Lambda=\bar G\cap\hat\Lambda,
\]
and $\bar\Lambda$ is a closed cocompact subgroup of $\bar G$.
\end{lemma}

\begin{proof}
The equality $\bar\Lambda=\bar G\cap\hat\Lambda$ follows directly from the derivative characterization of polynomial maps. As in \Cref{prop:hat G}, the map
\[
\bar G/\bar \Lambda
\longrightarrow
\Poly(\Z^\omega,G_{\bullet+1})/
\Poly(\Z^\omega,\Lambda_{\bullet+1}),
\qquad
(p,\gamma)\bar \Lambda
\longmapsto
p\Poly(\Z^\omega,\Lambda_{\bullet+1}),
\]
is a homeomorphism. The target is compact by \Cref{cor:shifted-quotients}, hence $\bar\Lambda$ is cocompact in $\bar G$.
\end{proof}

In the next two propositions and a lemma, we build the degree-$k$ translational structure into which we embed a polynomial orbit system. 

\begin{proposition}\label{prop:closed-orbits}
Let $x=h\hat \Lambda\in \hat G/\hat \Lambda$ and put
\[
K_x\coloneqq \mathrm{Stab}_{\bar G}(x).
\]
Then
\[
K_x=\bar G\cap h\hat \Lambda h^{-1},
\qquad
h^{-1}K_xh=\bar\Lambda.
\]
The subgroup $K_x$ is closed and cocompact in $\bar G$, the orbit $\bar Gx$ is compact and closed, and
\[
\bar G/K_x\longrightarrow\bar Gx,
\qquad
bK_x\longmapsto bx,
\]
is a homeomorphism.
\end{proposition}

\begin{proof}
An element $b\in\bar G$ belongs to $K_x$ if and only if $h^{-1}bh\in\hat\Lambda$, and hence
\[
K_x=\bar G\cap h\hat\Lambda h^{-1}.
\]
Since $\bar G$ is normal in $\hat G$, it follows that
\[
h^{-1}K_xh
=
\bar G\cap\hat\Lambda
=
\bar\Lambda.
\]
Thus conjugation by $h$ carries $\bar\Lambda$ isomorphically onto $K_x$. By \Cref{lem:bar Lambda}, $K_x$ is closed and cocompact in $\bar G$.

By normality of $\bar G$,
\[
\bar Gx
=
\bar Gh\hat\Lambda/\hat\Lambda
=
h\bar G\hat\Lambda/\hat\Lambda.
\]
The orbit map $b\mapsto bx$ induces a continuous bijection
\[
\bar G/K_x\longrightarrow\bar Gx.
\]
Since $\bar G/K_x$ is compact, this map is a homeomorphism. Hence $\bar Gx$ is compact, and therefore closed in $\hat G/\hat\Lambda$.
\end{proof}

\begin{lemma}\label{lem:adapted}
For every $x=h\hat\Lambda\in\hat G/\hat\Lambda$, the filtration $\bar G_\bullet$ is adapted to $K_x$.
\end{lemma}

\begin{proof}
By \Cref{prop:bar G}, each $\bar G_i$ is normal in $\hat G$. Hence
\[
h^{-1}(K_x\cap\bar G_i)h
=
\bar\Lambda\cap\bar G_i.
\]
For $i=0,1$, the quotient
\[
\bar G_i/(\bar G_i\cap\bar\Lambda)
=
\bar G/\bar\Lambda
\]
is compact by \Cref{lem:bar Lambda}. For $i\geq2$,
\[
\bar G_i/(\bar G_i\cap\bar\Lambda)
\simeq
\Poly(\Z^\omega,G_{\bullet+i})/
\Poly(\Z^\omega,\Lambda_{\bullet+i}),
\]
which is compact by \Cref{cor:shifted-quotients}. Conjugation by $h$ induces a homeomorphism
\[
\bar G_i/(\bar G_i\cap K_x)
\simeq
\bar G_i/(\bar G_i\cap\bar\Lambda).
\]
Thus $\bar G_i\cap K_x$ is cocompact in $\bar G_i$ for every $i\geq0$.

Moreover, the image of the compact space $\bar G_i/(\bar G_i\cap K_x)$ in $\bar G/K_x$ is compact and hence closed. Its inverse image under the quotient map $\bar G\to\bar G/K_x$ is $\bar G_iK_x$, so $\bar G_iK_x$ is closed in $\bar G$.
\end{proof}

\begin{proposition}\label{prop:translational-orbit}
Let $x\in\hat G/\hat\Lambda$. Then $\bar G/K_x$ can be equipped with the structure of a degree-$k$ translational $\Z^\omega$-system for the homomorphism
\[
\phi\colon\Z^\omega\to\bar G,
\qquad
\phi(\gamma)\coloneqq(1,\gamma).
\]
\end{proposition}

\begin{proof}
By \Cref{prop:bar G}, $\bar G_\bullet$ is a filtration of degree at most $k$. By \Cref{lem:adapted}, $K_x$ is adapted to $\bar G_\bullet$. 
\end{proof}

\Cref{prop:translational-orbit} gives a degree-$k$ translational structure on each compact homogeneous orbit $\bar Gx$. We now use this to embed polynomial orbit systems into translational systems.

\begin{theorem}\label{thm:polynomial-orbit-system}
Let $(G/\Lambda,G_\bullet)$ be a degree-$k$ filtered nilmanifold and let $p\in\Poly(\Z^\omega,(G/\Lambda)_\bullet)$. 
Then the associated polynomial orbit $\Z^\omega$-system $(Y_p,\sigma)$ is equivariantly homeomorphic to a closed invariant topological subsystem of a degree-$k$ translational $\Z^\omega$-system. 
\end{theorem}

\begin{proof}
Put
\[
x_p\coloneqq\xi^{-1}(p),
\qquad
K_p\coloneqq\mathrm{Stab}_{\bar G}(x_p).
\]
By \Cref{prop:translational-orbit}, $\bar G/K_p$ has the structure of a degree-$k$ translational $\Z^\omega$-system for the homomorphism
\[
\phi\colon\Z^\omega\to\bar G,
\qquad
\phi(\gamma)=(1,\gamma).
\]
By the equivariance in \Cref{prop:hat G},
\[
\xi^{-1}(\sigma^\gamma(p))
=
(1,\gamma)x_p.
\]
Therefore
\[
\xi^{-1}(Y_p)
=
\overline{\{(1,\gamma)x_p:\gamma\in\Z^\omega\}}
\subseteq
\bar Gx_p.
\]
By \Cref{prop:closed-orbits}, $\bar Gx_p\simeq\bar G/K_p$, and the displayed set is closed and invariant. Hence the restriction of $\xi^{-1}$ to $Y_p$, followed by this identification, is an equivariant homeomorphism from $(Y_p,\sigma)$ onto a closed invariant topological subsystem of $\bar G/K_p$.
\end{proof}

\section{Unique ergodicity in translational systems} \label{sec:unique ergodicity}

The goal of this section is to prove that orbit-closures in translational systems are uniquely ergodic.

\begin{theorem} \label{thm:equidistribution}
    Let $(G/\Lambda, \phi)$ be a translational $\Gamma$-system.
    Then for every $x \in G/\Lambda$, the orbit-closure $\mathcal{O}(x) = \overline{\{\phi(\gamma)x : \gamma \in \Gamma\}}$ is uniquely ergodic.
\end{theorem}

We remark that combining \Cref{thm:equidistribution} with \Cref{thm:polynomial-orbit-system} yields \Cref{thm:minimal finite order systems are u.e.}.

\begin{proof}[Proof of \Cref{thm:minimal finite order systems are u.e.}]
    Let $(X, T)$ be a minimal topological $\Gamma$-system of order $k$.
    Let $q \colon \Z^{\omega} \to \Gamma$ be a surjective group homomorphism, and let $S \colon \Z^{\omega} \acts X$ be the $\Z^{\omega}$ action $S^v = T^{q(v)}$ for $v \in \Z^{\omega}$.
    A measure on $X$ is $T$-invariant if and only if it is $S$-invariant, so it suffices to prove that $(X, S)$ is uniquely ergodic.

    Note that $(X, S)$ is a minimal $\Z^{\omega}$-system of order $k$, so by \cite[Theorem 1.5]{c_g-s_sz}, it is a factor of an inverse limit of minimal polynomial orbit $\Z^{\omega}$-systems.
    Unique ergodicity is preserved by factors and inverse limits, so we may assume that $(X, S)$ is itself a polynomial orbit $\Z^{\omega}$-system.
    By \Cref{thm:polynomial-orbit-system}, we may realize $(X,S)$ as a minimal component of a degree $k$ translational system and such minimal components are uniquely ergodic by \Cref{thm:equidistribution}.
\end{proof}

As described in the introduction, we prove \Cref{thm:equidistribution} by induction on degree, using an induction argument interweaving the unique ergodicity statement of \Cref{thm:equidistribution} with a Wiener--Wintner-type result to leverage the nilpotent structure of translational systems. 

This Wiener--Wintner type theorem needs weights that generalize (basic) nilsequences in the following sense: 

\begin{definition}[$\Gamma$-nilsequence]
    Let $\Gamma$ be a countable abelian group, and let $k \in \N$.
    A \emph{degree $k$ $\Gamma$-nilsequence} is a function $\psi \colon \Gamma \to \C$ of the form $\psi(\gamma) = F(\phi(\gamma)\Lambda)$ for a degree $k$ translational $\Gamma$-system $(G/\Lambda, \phi)$ and a continuous function $F \in C(G/\Lambda)$.
\end{definition}

\begin{remark}
    Note that the degree $k$ $\Gamma$-nilsequences are examples of $\Gamma$-sequences of order $k$ as defined in the introduction.
    In the case $\Gamma=\Z^{\omega}$, it follows from \cite[Theorem 1.5]{c_g-s_sz} and \Cref{thm:polynomial-orbit-system} that every $\Z^{\omega}$-sequence of order $k$ is a uniform limit of degree $k$ $\Z^{\omega}$-nilsequences.
    For arbitrary abelian groups $\Gamma$, the question of whether every $\Gamma$-sequence of order $k$ is approximable by degree $k$ $\Gamma$-nilsequences is intimately related to the problem of proving a structure theorem for $\Gamma$-systems of order $k$.
\end{remark}

The Wiener--Wintner type statement that we prove in conjunction with \Cref{thm:equidistribution} is the following:

\begin{theorem} \label{thm:WW-translational}
    Let $\Phi$ be a tempered F{\o}lner sequence in $\Gamma$.
    Suppose $(X, \mu, T)$ is a measure-preserving $\Gamma$-system.
    Then for every $f \in L^1(\mu)$, there exists a subset $X'_f \subseteq X$ with $\mu(X'_f) = 1$ such that if $x \in X'_f$ and $\psi \colon \Gamma \to \C$ is a $\Gamma$-nilsequence, then
    \begin{equation*}
        \frac{1}{|\Phi_N|} \sum_{\gamma \in \Phi_N} \psi(\gamma) f(T^{\gamma}x)
    \end{equation*}
    converges.
\end{theorem}

Using a similar argument to the one above by which we deduced \Cref{thm:minimal finite order systems are u.e.} from \Cref{thm:equidistribution}, we can also deduce \Cref{thm:WW} from \Cref{thm:WW-translational}.
The reliance on tempered F{\o}lner sequences creates a small technical hurdle, which is overcome by the following lemma.

\begin{lemma} \label{lem:Folner lifting}
	Let $\Gamma$ be a countable abelian group, and let $\Phi$ be a F{\o}lner sequence in $\Gamma$.
	Suppose $\tilde{\Gamma}$ is another countable abelian group and $q : \tilde{\Gamma} \to \Gamma$ is a surjective group homomorphism.
	Then there exists a F{\o}lner sequence $\tilde{\Phi}$ in $\tilde{\Gamma}$ such that $q(\tilde{\Phi}_N) = \Phi_N$ for every $N \in \N$, and $|q^{-1}(\{\gamma\}) \cap \tilde{\Phi}_N|$ is constant in $\gamma \in \Phi_N$ (depending only on $N$ and not on $\gamma$).
	In particular, for any bounded function $f : \Gamma \to \C$,
	\begin{equation*}
		\frac{1}{|\Phi_N|} \sum_{\gamma \in \Phi_N} f(\gamma) = \frac{1}{|\tilde{\Phi}_N|} \sum_{\tilde{\gamma} \in \tilde{\Phi}_N} (f \circ q)(\tilde{\gamma}).
	\end{equation*}
	Moreover, if $\Phi$ is tempered, then we can choose $\tilde{\Phi}$ to also be tempered.
\end{lemma}

\begin{proof}
	Let $K = \ker(q) \le \Gamma$.
	If $K$ is finite, we can take $\tilde{\Phi}_N = q^{-1}(\Phi_N)$.
	When $K$ is infinite, we must carefully choose a F{\o}lner sequence in the fibers instead.
	
	Let $t_1, t_2, \ldots$ be an enumeration of $\tilde{\Gamma}$.
	Let $s : \Gamma \to \tilde{\Gamma}$ be a section of $q$.
	Put $\Psi_1 = \{0\}$ and $\tilde{\Phi}_1 = s(\Phi_1)$.
	For $N \geq 2$, let
	\begin{equation*}
		V_N = \left\{ s(\gamma) + t_i - s(\gamma + q(t_i)) : i \leq N, \gamma \in \bigcup_{j=1}^N \Phi_j \right\} \subseteq K.
	\end{equation*}
	By amenability of $K$, we may find a set $\Psi_N \subseteq K$ such that
	\begin{equation*}
		\frac{\left| \left( \Psi_N + V_N \right) \triangle \Psi_N \right|}{|\Psi_N|} < 2^{-N},
	\end{equation*}
	Define $\tilde{\Phi}_N = s(\Phi_N) + \Psi_N$.
	By construction, $q(\tilde{\Phi}_N) = \Phi_N$, and each fiber has $|\Psi_N|$ elements.
	
	Let us check that $\tilde{\Phi}$ is a F{\o}lner sequence.
	Given $t \in \tilde{\Gamma}$, we have $t = t_n$ for some $n \in \N$.
	Now for $N > n$,
	\begin{align*}
		\tilde{\Phi}_N + t & = s(\Phi_N) + \Psi_N + t_n \\
		 & = \left\{ s(\gamma) + x + t_n : \gamma \in \Phi_N, x \in \Psi_N \right\} \\
		 & = \left\{ s(\gamma + q(t_n)) + \underbrace{\left( s(\gamma) + t_n - s(\gamma+q(t_n)) \right)}_{\in V_N} + x : \gamma \in \Phi_N, x \in \Psi_N \right\}.
	\end{align*}
	Therefore,
	\begin{equation*}
		\frac{\left| (\tilde{\Phi}_N + t) \cap \tilde{\Phi}_N \right|}{\left| \tilde{\Phi}_N \right|} \geq \underbrace{\frac{\left| (\Phi_N + q(t_n)) \cap \Phi_N \right|}{|\Phi_N|}}_{\to 1} \cdot \underbrace{\min_{v \in V_N} \frac{\left| (\Psi_N + v) \cap \Psi_N \right|}{|\Psi_N|}}_{> 1 - 2^{-N}} \to 1,
	\end{equation*}
	so $\tilde{\Phi}$ is a F{\o}lner sequence.
	
	Now suppose $\Phi$ is tempered, say with $|\Phi_N - \bigcup_{n=1}^{N-1} \Phi_n| \leq C|\Phi_N|$ for $N \in \N$.
	For each $N \in \N$, let
	\begin{equation*}
		W_N = V_N \cup \left\{ s(\gamma) - s(\gamma') - s(\gamma - \gamma') - y : \gamma \in \Phi_N, (\gamma', y) \in \bigcup_{n=1}^{N-1} (\Phi_n \times \Psi_n) \right\},
	\end{equation*}
	and define $\Psi_N$ and $\tilde{\Phi}_N$ using $W_N$ in place of $V_N$.
	Then $\tilde{\Phi}$ is still a F{\o}lner sequence by the argument above.
	Moreover, each element $u \in \tilde{\Phi}_N - \bigcup_{n=1}^{N-1} \tilde{\Phi}_n$ is of the form
	\begin{equation*}
		u = s(\gamma) + x - \left( s(\gamma') + y \right)
	\end{equation*}
	with $\gamma \in \Phi_N$, $x \in \Psi_N$, $\gamma' \in \Phi_n$, and $y \in \Psi_n$ for some $n \leq N-1$.
	Rearranging, we have
	\begin{equation*}
		u = s(\gamma - \gamma') + x + s(\gamma) - s(\gamma') - s(\gamma - \gamma') - y
		 \in s \left( \Phi_N - \bigcup_{n=1}^{N-1} \Phi_n \right) + \Psi_N + W_N,
	\end{equation*}
	so
	\begin{equation*}
		\left| \tilde{\Phi}_N - \bigcup_{n=1}^{N-1} \tilde{\Phi}_n \right| \leq \underbrace{\left| \Phi_N - \bigcup_{n=1}^{N-1} \Phi_n \right|}_{\leq C |\Phi_N|} \cdot \underbrace{\left| \Psi_N + W_N \right|}_{\leq \left( 1 + 2^{-N} \right) |\Psi_N|} \leq 2C \left| \tilde{\Phi}_N \right|.
	\end{equation*}
	Thus, $\tilde{\Phi}$ is tempered.
\end{proof}

\begin{proof}[Proof that \Cref{thm:WW-translational} implies \Cref{thm:WW}]
    Suppose \Cref{thm:WW-translational} holds for the group $\Z^{\omega}$.
    Let $\Gamma$ be an arbitrary countable abelian group, and let $\psi : \Gamma \to \C$ be a $\Gamma$-sequence of order $k$.
    Let $(X, \mu, T)$ be a measure-preserving $\Gamma$-system, $f \in L^1(\mu)$, and $\Phi$ a tempered F{\o}lner sequence in $\Gamma$.

    There exists a surjective homomorphism $q : \Z^{\omega} \to \Gamma$.
    The $\Z^{\omega}$-sequence $\tilde{\psi} = \psi \circ q$ is of order $k$.
    Hence, arguing as in the proof of \Cref{thm:minimal finite order systems are u.e.} above, $\tilde{\psi}$ can be uniformly approximated by degree $k$ $\Z^{\omega}$-nilsequences, so we may assume that $\tilde{\psi}$ is itself a $\Z^{\omega}$-nilsequence.
    Applying \Cref{lem:Folner lifting}, we can find a tempered F{\o}lner sequence $\tilde{\Phi}$ in $\Z^{\omega}$ satisfying
    \begin{equation*}
        \frac{1}{|\Phi_N|} \sum_{\gamma \in \Phi_N} \psi(\gamma) f(T^{\gamma}x) = \frac{1}{|\tilde{\Phi}_N|} \sum_{v \in \tilde{\Phi}_N} \tilde{\psi}(v) f(\tilde{T}^vx),
    \end{equation*}
    where $\tilde{T}^v = T^{q(v)}$.
    The average on the right-hand side converges almost everywhere by \Cref{thm:WW-translational}, so the average on the left-hand side also converges almost everywhere, proving \Cref{thm:WW}.
\end{proof}

For each $k \in \N$, we will work with the statements:

\begin{named}{Equi($k$)}{} \label{Equi}
    \Cref{thm:equidistribution} holds for degree $k$ translational systems.
\end{named}

\begin{named}{WW($k$)}{} \label{WW}
    \Cref{thm:WW-translational} holds for degree $k$ nilsequences.
\end{named}

The induction strategy is to prove the implications $$\textup{Equi}(k) \implies \textup{WW}(k) \implies \textup{Equi}(k+1)$$ for every $k \in \N$.

Let us remark that WW($1$) was previously established in \cite[p. 121]{ornstein-weiss} and \cite[Corollary 4.1]{zorin_ww} and used to prove Equi($2$) in \cite[\S 3]{AJ}.
It does not provide any advantage to start our induction from the statement Equi($2$) rather than from Equi($1$), so we will in the course of our proof also arrive at new (though largely similar to the previous) proofs of these statements.

A useful observation is that, by the same argument used to deduce \Cref{thm:minimal finite order systems are u.e.} from \Cref{thm:equidistribution} above, the statement \ref{Equi} implies that minimal topological $\Gamma$-systems of order $k$ are uniquely ergodic.

\bigskip

The next three subsections are dedicated to technical results that will come into play in the proof of \Cref{thm:equidistribution}, which then appears in the last two subsections.

\subsection{Orbit-closures in degree $k$ translational systems are systems of order $k$}

\begin{proposition} \label{prop:orbits order k}
    Let $(G/\Lambda, \phi)$ be a degree $k$ translational $\Gamma$-system.
	If $Y \subseteq G/\Lambda$ is a minimal set, then $(Y, \phi)$ is a topological system of order $k$. 
\end{proposition}

\begin{proof}
    We want to show that the $k$-regionally proximal relation on $Y$ is trivial.
    By \cite[Proposition 2.3.1]{candela0}, $X = G/\Lambda$ equipped with the cubes $C^n(X) = \{\pi_{\Lambda} \circ c : c \in C^n(G_{\bullet})\}$ is a $k$-step nilspace, i.e. every $(k+1)$-corner has unique completion.
    The collection of dynamical $(k+1)$-cubes
    \begin{equation*}
        Q^{\bracket{k+1}}(Y) = \overline{\left\{ \left( \phi(\bm{\omega} \cdot \bm{\gamma}) y \right)_{\bm{\omega} \in \bracket{k+1}} : y \in Y, \bm{\gamma} = (\gamma_1,\ldots,\gamma_{k+1}) \in \Gamma^{k+1} \right\}}
    \end{equation*}
    is a subset of the $(k+1)$-cubes $C^{k+1}(X)$, so the dynamical cubes also have unique corner completion.
    That is, the $k$-regionally proximal relation is trivial.
\end{proof}

\begin{corollary} \label{cor:orbits measurably order k}
    Let $(G/\Lambda, \phi)$ be a degree $k$ translational $\Gamma$-system.
	If $Y \subseteq G/\Lambda$ is a minimal set and $\nu$ is an ergodic invariant measure supported on $Y$, then $(Y, \nu, \phi)$ is a measure-preserving system of order $k$. 
\end{corollary}

\begin{proof}
    This follows from \Cref{prop:orbits order k} and \Cref{prop: topological order k implies measurable order k}.
\end{proof}

\subsection{Degree reduction of $\Gamma$-nilsequences}

\begin{definition}[Vertical characters]
    Let $G/\Lambda$ be a degree $k$ filtered translational space.
    A continuous function $f \colon G/\Lambda \to S^1$ is a \emph{vertical character} if there exists a character $\chi \in \hat{G}_k$ such that $f(ux) = \chi(u)f(x)$ for every $u \in G_k$ and $x \in G/\Lambda$.
\end{definition}

\begin{proposition} \label{prop:characters dense}
    Let $G/\Lambda$ be a finite degree filtered translational space.
    Then linear combinations of vertical characters are uniformly dense in $C(G/\Lambda)$.
\end{proposition}

\begin{proof}
This was established in \cite[Proposition 4.2]{Holikov-Shalom} for nilmanifolds, but the same proof holds in our case. We give the full proof for the sake of completeness. Let $\mathcal{A}$ denote the closed algebra of all continuous vertical characters. By the Stone--Weierstra\ss~theorem, it suffices to show that $\mathcal{A}$ separate points. Let $G/\Lambda$ be a translational system of degree $k$ for some $k\geq 1$, and let $x,y\in G/\Lambda$ be two distinct points. If $\pi(x)\not= \pi(y)$, where $\pi\colon G/\Lambda\rightarrow G/G_k\Lambda$ is the natural quotient map, then by Uryshon's lemma we can find a continuous $f\colon G/G_k\Lambda\rightarrow S^1$ with $f(\pi(x))=1,f(\pi(y))=-1$. Since $f\circ\pi\in \mathcal{A}$ (a vertical character with trivial frequency), we deduce that $\mathcal{A}$ separates $x$ and $y$ in this case. Otherwise, $\pi(x)=\pi(y)$ and so there exists $u\in G_k$ so that $x=uy$. Since $\Lambda_k:=\Lambda\cap G_k$ acts trivially on $G/\Lambda$ and $x \ne y$, the coset $u\Lambda_k$ is non-trivial. In particular, since the characters on $G_k/\Lambda_k$ separate points, there is some non-trivial $\chi\in \widehat{G_k/\Lambda_k}$ satisfying that $\chi(u\Lambda_k)\not = 1$. Again by Uryshon's lemma we can find a continuous bump function $f\colon G/\Lambda\rightarrow S^1$ in a small neighborhood around $x$ (depending on $\chi$). Next, define
\[
\phi\colon G/\Lambda\rightarrow \C \qquad \phi(x):=\int_{G_k/\Lambda_k} f(tx)\overline{\chi(t)}~dt
\]
where $dt$ is the Haar measure. The continuity of $f,\chi$ implies the continuity of the convolution $\phi$ and $\phi(x)\not=1$ by the choice of $f$. For every $v\in G_k/\Lambda_k$, we then have $$\phi(vx) = \int_{G_k/\Lambda_k} f(tvx)\overline{\chi(t)}dt =\int_{G_k/\Lambda_k} f(tx)\overline{\chi(tv^{-1})}dt = \chi(v)\phi(x),$$ whence $\phi\in A$ and $\phi(y)=\phi(ux)=\chi(u)\phi(x)\not=\phi(x).$ This completes the proof.
\end{proof}

The rest of this section is devoted to the proof of the following theorem which asserts that the nilsequence associated with the discrete derivative of a vertical character is a nilsequence of a smaller degree.
\begin{theorem}\label{thm:nilsequence}
    Let $k \in \N$, let $(G/\Lambda,\phi)$ be a degree $k+1$ translational $\Z^{\omega}$-system, and let $x_0 \in G/\Lambda$.
    Let $f\in C(X)$ be a vertical character.
    Then for every $g\in G$, the function
    \begin{equation*}
        \psi_g(\gamma):= f(\phi(\gamma)gx_0)\cdot\overline{f(\phi(\gamma)x_0)}.
    \end{equation*}
    is a $\Z^\omega$-nilsequence of degree $k$.
\end{theorem}
The proof requires building a translational system where one can realize $\psi_g(\gamma)$ as a nilsequence. Our first goal is to construct the translational system. Inspired by \cite{gt1} we define $$G^{\square}=G\times_{G_2}G = \{(g_1,g_2)\in G^2 : g_1^{-1}g_2\in G_2\}.$$
It is obvious that $G^{\square}\leq G\times G$ is closed, and therefore a Polish group. 
Let $D = \diag(G) = \{(g,g) \in G \times G : g \in G\} \subseteq G^{\square}$, and let $N = \{1\} \times G_2 \subseteq G^{\square}$.
\begin{lemma}
    $G^{\square}=D\cdot N$.
    Moreover, for every $1\leq i \leq k$, let $G_\bullet=(G_i)_{i\geq 0}$ denote the filtration on $G$ and let $$G_i^{\square}\coloneqq D_i\cdot N_i$$ where $D_i=\diag(G_i)$ and $N_i=\{1\}\times{G_{i+1}}.$ Then, $G_\bullet^{\square}=(G_i^{\square})_{i\geq 0}$ is a filtration on $G^{\square}.$
\end{lemma}
\begin{proof}
    Since $D,N\subseteq G^{\square}$, we have $D\cdot N\subseteq G^{\square}$. In the other direction, let $(g_1,g_2)\in G^{\square}$ and write $(g_1,g_2)=(g_1,g_1)\cdot (1,g_1^{-1}g_2)\subseteq D\cdot N.$

Now, for all $i$ we prove that $D_i\cdot N_i$ is a group. To this end, let $(g,g),(h,h)\in D_i$ and $(1,u),(1,v)\in N_i$ and compute 
\begin{multline*}
    (g,g)(1,u)(h,h)(1,v) = (gh,guhv)=(gh,ghu[u,h]v)\\=(gh,gh)(1,u[u,h]v)\in D_i\cdot N_i
\end{multline*}
It is clearly closed, so it suffices to show that $[G^\square_i,G_j^\square] \subseteq G_{i+j}^\square.$ Let $(g,g)\in D_i$, $(1,u)\in N_i$ and $(h,h)\in D_j$, $(1,v)\in N_j$. Then,
$$[(g,gu),(h,hv)]=([g,h],[gu,hv]) = ([g,h],[g,h]) \cdot (1, [g,h]^{-1}[gu,hv]).$$
Since $([g,h],[g,h])\in D_{i+j}$, it suffices to show that $[g,h]^{-1}[gu,hv]\in G_{i+j+1}.$ Modulo $G_{i+j+1}$ we have that $u$ commutes with $h,v$ and $v$ commutes with $g,u$ and both $u$ and $v$ commute with $[g,h]$. Therefore,
$$[gu,hv]=u^{-1}g^{-1}v^{-1}h^{-1}guhv \equiv u^{-1}v^{-1}[g,h]uv  \equiv [g,h][u,v] \pmod{G_{i+j+1}}.$$
We conclude that $[g,h]^{-1}[gu,hv] \equiv [u,v] \equiv 0 \pmod{G_{i+j+1}}$, hence the claim follows.
\end{proof}
Recall that $\Gamma$ acts on $G/\Lambda$ via a homomorphism $\phi\colon\Gamma\rightarrow G$. This extends to a homomorphism $\phi^{\square}\colon\Gamma\rightarrow G^{\square}$ by taking $\phi^{\square}(\gamma)=(\phi(\gamma),\phi(\gamma))\in G^{\square}.$\\

Next we construct the lattice. Let $x_0 \in X = G/\Lambda$, and let $g\in G$ be arbitrary as in the theorem. Set $z_g\coloneqq (x_0,gx_0)\in X\times X$. Let $\Lambda_g \coloneqq \mathrm{stab}_{G^\square}(z_g).$
\begin{lemma}\label{orbitistranslational}
The subset $Y_g=G^{\square}z_g = \{hz_g : h\in G^{\square}\}$ is isomorphic to the translational system $G^\square/\Lambda_g$ as a topological dynamical system.
\end{lemma}
\begin{proof}
    Since $G^{\square}\leq G\times G$ is the kernel of the homomorphism $(g_1,g_2)\mapsto g_1^{-1}g_2G_2$, it is closed (since $G/G_2$ is Hausdorff) and normal. Let $q\colon G\times G\rightarrow X\times X$ denote the quotient map (by $\Lambda\times \Lambda)$. We analyze $q^{-1}(Y_g).$

    Observe that $z_g = q(1,g)$ and recall that every element of $G^\square$ takes the form $(h,hu)$ with $h\in G$ and $u\in G_2$. Since $q^{-1}\circ q(S) = S\cdot (\Lambda\times \Lambda)$ for any subset $S\subseteq G$, we have
    $$q^{-1}(Y_g) = G^\square(1,g)(\Lambda\times \Lambda) = \{(h\lambda_1,hug\lambda_2):h\in G,\, u\in G_2\, \lambda_1,\lambda_2\in \Lambda\}.$$

    By definition, the action of $G^{\square}$ on $Y_g$ is transitive. In order to apply Effros's theorem we must first ensure that $Y_g$ is compact. For this purpose, observe that the map $\sigma(a,b)=(a,ab)$ is a homeomorphism of $G\times G$ and its inverse is given by $(a,b)\mapsto (a,a^{-1}b)$. Therefore, $$\sigma^{-1}(h\lambda_1,hug\lambda_2) = (h\lambda_1, \lambda_1^{-1}ug\lambda_2),$$ which implies that $\sigma^{-1}(q^{-1}(Y_g)) = G\times \Lambda G_2g\Lambda =G\times  gG_2\Lambda$ by normality of $G_2 \leq G$.
    By assumption $G_2\Lambda$ is closed and therefore so is $gG_2\Lambda$. Since $q^{-1}(Y_g)$ is closed, the complement $(G\times G)\backslash q^{-1}(Y_g)$ is open. Since the quotient is an open map, $(X\times X)\backslash Y_g$ is open from which we deduce that $Y_g$ is closed, and therefore compact (since $X\times X$ is Hausdorff).

    In particular, $Y_g$ is a Polish group and so by Effros's theorem $Y_g\cong G^\square/\Lambda_g.$ Hence, $Y_g$ is a translational system.
\end{proof}

The proof of \Cref{thm:nilsequence} is now simple.

\begin{proof}[Proof of \Cref{thm:nilsequence}]
    Let $f$ be a vertical character with frequency $\chi$. Define the function $F\colon X\times X\rightarrow S^1$ by setting $F(x,y)=f(y)\cdot\overline{f(x)}$. The continuity of $f$ implies that $F$ is continuous and moreover it is invariant to the action of $D_{k+1} = G^\square_{k+1}$. Restricting $F$ to $Y_g\cong G^\square/\Lambda_g$ we see that $F$ factors through the quotient $Y'_g=G^\square/G^\square_k\Lambda_g$ that is a $k$-step translational system. Hence, $F(\phi^\square (\gamma)z_g) = F(\phi(\gamma) x_0,\phi(\gamma) gx_0) =f(\phi(\gamma) gx_0)\cdot \overline{f(\phi(\gamma) x_0)}$ is a $k$-step nilsequence.
\end{proof}

\subsection{$\textup{Equi}(k) \implies \textup{WW}(k)$}

We now prove the first part of the induction step, upgrading the unique ergodicity statement for orbits in degree $k$ translational system to a Wiener--Wintner theorem for degree $k$ nilsequence weights.
Our approach follows the well-trod method of \emph{characteristic factors}.
We prove that an ergodic average weighted by a degree $k$ nilsequence is ''controlled'' by the order $k$ Host--Kra factor.
This essentially reduces the task of proving WW($k$) to the special case when the measure-preserving system is of order $k$.
We then use recent results from nilspace theory to produce a topological model that allows us deduce convergence of ergodic averages using Equi($k$).

\bigskip

The next proposition says that the order $k$ Host--Kra factor is characteristic for ergodic averages with degree $k$ nilsequence weights.

\begin{proposition} \label{prop:WW characteristic factor}
    Let $(X, \mu, T)$ be an ergodic measure-preserving system.
    Let $\Phi$ be a F{\o}lner sequence in $\Gamma$, and suppose $x \in \gen(\mu,\Phi)$.
    Let $k \in \N$, and assume that \ref{Equi} holds.
    Then for every $f \in C(X)$ and every degree $k$ $\Gamma$-nilsequence $\psi \colon \Gamma \to \C$,
    \begin{equation*}
        \limsup_{N \to \infty} \left| \frac{1}{|\Phi_N|} \sum_{\gamma \in \Phi_N} \psi(\gamma) f(T^{\gamma}x) \right| \leq \|\psi\|_{\infty} \cdot \left\| \E[f \mid Z_k] \right\|_{L^1(\mu)}
    \end{equation*}
\end{proposition}

\begin{proof}
    We follow the proof strategy from \cite[Lemma 6.11]{sumsets}.

    Let $(G/\Lambda, \phi)$ be a degree $k$ translational $\Gamma$-system and $F \in C(G/\Lambda)$ such that $\psi(\gamma) = F(\phi(\gamma)\Lambda)$.
    We may assume $\|F\|_{\sup} = \|\psi\|_{\infty}$.
    Let $Y = \overline{\{\phi(\gamma)\Lambda : \gamma \in \Gamma\}} \subseteq G/\Lambda$.
    By \ref{Equi}, the system $(Y, \phi)$ is uniquely ergodic.
    Let $\nu$ be the unique invariant measure on $Y$.

    Let $(N_k)_{k \in \N}$ be a subsequence such that
    \begin{equation*}
        \limsup_{N \to \infty} \left| \frac{1}{|\Phi_N|} \sum_{\gamma \in \Phi_N} \psi(\gamma) f(T^{\gamma}x) \right| = \lim_{k \to \infty} \left| \frac{1}{|\Phi_{N_k}|} \sum_{\gamma \in \Phi_{N_k}} \psi(\gamma) f(T^{\gamma}x) \right|.
    \end{equation*}
    By passing to a further subsequence, we may assume that the sequence of measures
    \begin{equation*}
        \frac{1}{|\Phi_{N_k}|} \sum_{\gamma \in \Phi_{N_k}} \delta_{T^{\gamma}x} \times \delta_{\phi(\gamma)\Lambda}
    \end{equation*}
    converges in the weak* topology to some measure $\rho$ on $X \times Y$.
    Since $x$ is generic for $\mu$ along $\Phi$ and $\Lambda \in Y$ is generic for $\nu$ (by unique ergodicity of $Y$), the measure $\rho$ is a joining of $\mu$ and $\nu$.
    We will use this property to check that
    \begin{equation} \label{eq:joining project to factor}
        \int_{X \times Y} u \otimes v~d\rho = \int_{X \times Y} \E[u \mid Z_k] \otimes v~d\rho
    \end{equation}
    for every $u \in L^2(\mu)$ and $v \in L^2(\nu)$.

    To prove \eqref{eq:joining project to factor}, we begin by making several reductions.
    First, by linearity of the expression, we may assume that $u \in L^2(\mu)$ satisfies $\E[u \mid Z_k] = 0$.
    Next, by taking the ergodic decomposition of $\rho$, we may assume $\rho$ is an ergodic joining of $\mu$ and $\nu$.
    (The fact that the ergodic components of $\rho$ are themselves joinings is shown, for example, in \cite[Proposition 1]{delarue_joinings}.)

    Since $\rho$ is a joining of $\mu$ and $\nu$, we have
    \begin{equation*}
        \|u \otimes 1\|_{U^{k+1}(X \times Y, \rho, T \times \phi)} = \|u\|_{U^{k+1}(X,\mu,T)} = 0.
    \end{equation*}
    On the other hand, the function $1 \otimes v$ is measurable with respect to the factor $(Y, \nu, \phi)$, which is of order $k$ by \Cref{cor:orbits measurably order k}.
    Hence, $1 \otimes \nu$ is $Z_k(X \times Y, \rho, T \times \phi)$-measurable.
    Therefore,
    \begin{equation*}
        \int_{X \times Y} u \otimes v~d\rho = \int_{X \times Y} (u \otimes 1)(1 \otimes v)~d\rho = 0.
    \end{equation*}

    To finish, using \eqref{eq:joining project to factor} and H\"{o}lder's inequality, we have
    \begin{align*}
        \limsup_{N \to \infty} \left| \frac{1}{|\Phi_N|} \sum_{\gamma \in \Phi_N} \psi(\gamma) f(T^{\gamma}x) \right| & = \left| \int_{X \times Y} f \otimes F~d\rho \right| \\
        & = \left| \int_{X \times Y} \E[f \mid Z_k] \otimes F~d\rho \right| \\
        & \leq \|F\|_{\sup} \cdot \left\| \E[f \mid Z_k] \right\|_{L^1(\mu)} \\
        & = \|\psi\|_{\infty} \cdot \left\| \E[f \mid Z_k] \right\|_{L^1(\mu)}.
    \end{align*}
\end{proof}

\begin{proposition}
    For every $k \in \N$, $\textup{Equi}(k) \implies \textup{WW}(k)$.
\end{proposition}

\begin{proof}
    Assume \ref{Equi} holds.
    Let $(X, \mu, T)$ be a measure-preserving $\Gamma$-system.
    We will assume that $X$ is a compact metric space and $T$ acts by homeomorphisms $T^{\gamma} \colon X \to X$.
    Using the ergodic decomposition, we may also assume that $\mu$ is an ergodic measure.
    Let $\Phi$ be a tempered F{\o}lner sequence in $\Gamma$.
    By the maximal ergodic theorem for amenable groups \cite[Theorem 3.2]{lindenstrauss}, it suffices to verify \ref{WW} for functions $f \in C(X)$.

    We will pass to an extension of $(X, \mu, T)$ in order to realize the order $k$ Host--Kra factor as a topological factor.
    By \cite[Theorem 5.11]{candela-szegedy-cubic} and \cite[Proposition 5.4]{CGSSextensions}, there exists a minimal, distal, uniquely ergodic system $(Z_k, \mu_k, T_k)$ that is isomorphic to the order $k$ Host--Kra factor of $(X, \mu, T)$.
    Let $\tilde{X} = X \times Z_k$, and let $\tilde{\mu}$ be the relatively independent joining of $\mu$ and $\mu_k$ over $Z_k$.
    That is,
    \begin{equation*}
        \int_{X \times Z_k} h_1 \otimes h_2~d\tilde{\mu} = \int_X h_1 \cdot (h_2 \circ \pi_k)~d\mu
    \end{equation*}
    for $h_1 \in C(X)$ and $h_2 \in C(Z_k)$, where $\pi_k$ is the measurable factor map $\pi_k \colon X \to Z_k$.
    Define an action $\tilde{T}$ by $\tilde{T}^{\gamma} = T^{\gamma} \times T_k^{\gamma} \colon \tilde{X} \to \tilde{X}$.
    Then $(\tilde{X}, \tilde{\mu}, \tilde{T})$ is isomorphic as a measure-preserving system to $(X, \mu, T)$, with the coordinate projection map $p \colon \tilde{X} \to X$ giving an isomorphism.
    Moreover, the other projection $q \colon \tilde{X} \to Z_k$ is a continuous factor map.

    Since $\Phi$ is tempered and $\tilde{\mu}$ is an ergodic measure, the set $\gen(\tilde{\mu},\Phi) \subseteq \tilde{X}$ of generic points for $\tilde{\mu}$ along $\Phi$ is a set of full measure by Lindenstrauss's pointwise ergodic theorem \cite[Theorem 1.2]{lindenstrauss}.
    
    Let $\tilde{f} = f \circ p \in C(\tilde{X})$.
    Approximating $\E[\tilde{f} \mid Z_k]$ by continuous functions in $L^1$, we have a decomposition
    \begin{equation*}
        \tilde{f} = \tilde{f}_{\textup{str}} \circ q + \tilde{f}_{\textup{unf}}
    \end{equation*}
    with $\tilde{f}_{\textup{str}} \in C(Z_k)$, $\tilde{f}_{\textup{unf}} \in C(\tilde{X})$, and $\left\|\tilde{f}_{\textup{str}} - \E[f \mid Z_k] \right\|_{L^1(\mu_k)} < \eps$ for arbitrary $\eps > 0$.
    (The continuity of the factor map $q$ is used to ensure continuity of $\tilde{f}_{\textup{unf}}$.)

    Now let $\psi \colon \Gamma \to \C$ be a degree $k$ $\Gamma$-nilsequence.
    By \Cref{prop:WW characteristic factor}, we have
    \begin{equation} \label{eq:uniform part}
        \limsup_{N \to \infty} \left| \frac{1}{|\Phi_N|} \sum_{\gamma \in \Phi_N} \psi(\gamma)  \tilde{f}_{\textup{unf}}(\tilde{T}^{\gamma}x) \right| \leq \eps \cdot \|\psi\|_{\infty}
    \end{equation}
    for every $x \in \gen(\tilde{\mu},\Phi)$.

    We now analyze the average
    \begin{equation*}
        \frac{1}{|\Phi_N|} \sum_{\gamma \in \Phi_N} \psi(\gamma) \tilde{f}_{\textup{str}}(T_k^{\gamma}z)
    \end{equation*}
    for $z \in Z_k$.
    Let $(G/\Lambda, \phi)$ be a degree $k$ translational $\Gamma$-system and $F \in C(G/\Lambda)$ such that $\psi(\gamma) = F(\phi(\gamma)\Lambda)$.
    We may assume $\|F\|_{\sup} = \|\psi\|_{\infty}$.
    By \ref{Equi}, the orbit-closure $Y_z = \overline{\{(\phi(\gamma)\Lambda, T_k^{\gamma}z) : \gamma \in \Gamma\}} \subseteq G/\Lambda \times Z_k$ is uniquely ergodic, say with unique invariant measure $\nu_z$.
    Therefore,
    \begin{equation*}
        \lim_{N \to \infty} \frac{1}{|\Phi_N|} \sum_{\gamma \in \Phi_N} \psi(\gamma) \tilde{f}_{\textup{str}}(T_k^{\gamma}z) = \int_{G/\Lambda \times Z_k} F \otimes \tilde{f}_{\textup{str}}~d\nu_z.
    \end{equation*}
    Since $Z_k$ is uniquely ergodic, the measure $\nu_z$ projects to $\mu_k$, so we have the estimate
    \begin{multline} \label{eq:structured part}
        \limsup_{N \to \infty} \left| \frac{1}{|\Phi_N|} \sum_{\gamma \in \Phi_N} \psi(\gamma) \tilde{f}_{\textup{str}}(T_k^{\gamma}z) - \int_{G/\Lambda \times Z_k} F \otimes \E[f \mid Z_k]~d\nu_z \right| \\
         = \left| \int_{G/\Lambda \times Z_k} F \otimes (\tilde{f}_{\textup{str}} - \E[f \mid Z_k])~d\nu_z \right| \leq \eps \cdot \|F\|_{\sup} = \eps \cdot \|\psi\|_{\infty}.
    \end{multline}
    
    Combining \eqref{eq:uniform part} and \eqref{eq:structured part}, if $x \in \gen(\tilde{\mu},\Phi)$, then
    \begin{align*}
        \limsup_{N \to \infty} & \left| \frac{1}{|\Phi_N|} \sum_{\gamma \in \Phi_N} \psi(\gamma)  \tilde{f}(\tilde{T}^{\gamma}x) - \int_{G/\Lambda \times Z_k} F \otimes \E[f \mid Z_k]~d\nu_{\pi_k(x)} \right| \\
         & \leq \limsup_{N \to \infty} \left| \frac{1}{|\Phi_N|} \sum_{\gamma \in \Phi_N} \psi(\gamma)  \tilde{f}_{\textup{str}}(T_k^{\gamma} \pi_k(x)) - \int_{G/\Lambda \times Z_k} F \otimes \E[f \mid Z_k]~d\nu_{\pi_k(x)} \right| \\
         & \qquad \qquad + \limsup_{N \to \infty} \left| \frac{1}{|\Phi_N|} \sum_{\gamma \in \Phi_N} \psi(\gamma)  \tilde{f}_{\textup{unf}}(T^{\gamma}x) \right| \\
         & \leq 2 \eps \cdot \|\psi\|_{\infty}.
    \end{align*}
    But $\eps > 0$ was arbitrary and $\psi$ is bounded, so we conclude that
    \begin{equation*}
        \frac{1}{|\Phi_N|} \sum_{\gamma \in \Phi_N} \psi(\gamma) \tilde{f}(\tilde{T}^{\gamma}x)
    \end{equation*}
    converges to
    \begin{equation*}
        \int_{G/\Lambda \times Z_k} F \otimes \E[f \mid Z_k]~d\nu_{\pi_k(x)}
    \end{equation*}
    for $x \in \gen(\tilde{\mu},\Phi)$.

    To finish, by Fubini's theorem, there exists a set $X' \subseteq X$ with $\mu(X') = 1$ such that for every $x \in X'$, there exists $\tilde{x} \in \gen(\tilde{\mu},\Phi)$ with $p(\tilde{x}) = x$.
    Hence, for $x \in X'$,
    \begin{equation*}
        \frac{1}{|\Phi_N|} \sum_{\gamma \in \Phi_N} \psi(\gamma) f(T^{\gamma}x) = \frac{1}{|\Phi_N|} \sum_{\gamma \in \Phi_N} \psi(\gamma) \tilde{f}(\tilde{T}^{\gamma}\tilde{x})
    \end{equation*}
    converges.
\end{proof}

\subsection{$\textup{WW}(k) \implies \textup{Equi}(k+1)$}

Given a translational $\Gamma$-system $(G/\Lambda, \phi)$ and a point $x \in G/\Lambda$, the orbit-closure $\mathcal{O}(x) = \overline{\{\phi(\gamma)x : \gamma \in \Gamma\}}$ is minimal.
This follows from the fact that translational systems are distal (see \cite[\S 2.4]{AJ}).
We can therefore utilize the following criterion for unique ergodicity.

\begin{lemma} \label{lem:oxtoby}
	Fix a F{\o}lner sequence $(\Phi_N)$ in $\Gamma$.
	A minimal topological dynamical $\Gamma$-system $(X,T)$ is uniquely ergodic if and only if
	\begin{equation*}
		\frac{1}{|\Phi_N|} \sum_{\gamma \in \Phi_N} f \left( T^{\gamma} x \right)
	\end{equation*}
	converges for every $f \in C(X)$ and every $x \in X$.
\end{lemma}

\begin{proof}
    This was established by Oxtoby \cite[Proposition 5.4]{oxtoby_ergodic} for $\Z$-systems and generalized to $\Gamma$-systems in \cite[Lemma 3.3]{AJ}.
\end{proof}

\begin{proposition}
    For every $k \in \N$, $\textup{WW}(k) \implies \textup{Equi}(k+1)$.
\end{proposition}

\begin{proof}
    Let $(G/\Lambda, \phi)$ be a degree $k+1$ translational $\Gamma$-system.
    The orbit-closures $\mathcal{O}(x)$, $x \in G/\Lambda$, are minimal systems, so it suffices by \Cref{lem:oxtoby} to find a F{\o}lner sequence $\Phi$ such that the ergodic averages
    \begin{equation} \label{eq:ergodic avg}
        \frac{1}{|\Phi_N|} \sum_{\gamma \in \Phi_N} f(T^{\gamma}x)
    \end{equation}
    converge for every $f \in C(G/\Lambda)$ and every $x \in G/\Lambda$.
    By \Cref{prop:characters dense}, it further suffices to prove convergence of \eqref{eq:ergodic avg} for vertical characters $f$.
    
    Fix a tempered F{\o}lner sequence $\Phi$ in $\Gamma$.
    Then by \ref{WW}, there exists $x_0 \in G/\Lambda$ such that
    \begin{equation*}
        \frac{1}{|\Phi_N|} \sum_{\gamma \in \Phi_N} \psi(\gamma) f(\phi(\gamma)x_0)
    \end{equation*}
    converges for every degree $k$ $\Gamma$-nilsequence $\psi$ and every $f \in C(G/\Lambda)$.
    Now, for arbitrary $x \in G/\Lambda$, we may write $x = gx_0$ for some $g \in G$.
    Given a vertical character $f \in C(G/\Lambda)$, we may then apply \Cref{thm:nilsequence} to write
    \begin{equation*}
        f(\phi(\gamma)x) = \psi(\gamma) f(\phi(\gamma)x_0)
    \end{equation*}
    for some degree $k$ $\Gamma$-nilsequence $\psi$ (which depends on the point $x$).
    Therefore,
    \begin{equation*}
        \frac{1}{|\Phi_N|} \sum_{\gamma \in \Phi_N} f(\phi(\gamma)x)
    \end{equation*}
    converges for every $x \in G/\Lambda$ and every vertical character $f \in C(G/\Lambda)$.
    This completes the proof.
\end{proof}

\section{Dynamical reformulations of the combinatorial results} \label{sec:dynamical reformulation}

In this section, we recast the main combinatorial results of the paper (\Cref{thm-main} and \Cref{thm-main1}) in a dynamical context.

\subsection{Erd\H{o}s progressions and Erd\H{o}s cubes}

Infinite sumset configurations in subsets of a group $\Gamma$ can be encoded by certain dynamical configurations in actions of $\Gamma$ by homeomorphisms on a compact metric space.
We begin by defining these dynamical objects.

\begin{definition}[Erd\H{o}s progressions]
    Let $(X, T)$ be a topological dynamical $\Gamma$-system, and let $k \in \N$.
    A tuple $(x_0, x_1, \dots, x_k) \in X^k$ is a \emph{($(k+1)$-term) Erd\H{o}s progression} if there exists a sequence $(\gamma_n)_{n \in \N}$ of distinct elements of $\Gamma$ such that
    \begin{equation*}
        \lim_{n \to \infty} \left( T^{\gamma_n} x_0, \dots, T^{\gamma_n} x_{k-1} \right) = (x_1, \dots, x_k).
    \end{equation*}
\end{definition}

\begin{definition}[Erd\H{o}s cubes] \label{def: erdos cube}
    Let $(X, T)$ be a topological dynamical $\Gamma$-system, and let $k \in \N$.
    A tuple $(x_{\bm{\omega}})_{\bm{\omega} \in \bracket{k}} \in X^{\bracket{k}}$ is a \emph{$k$-dimensional Erd\H{o}s cube} if there exists $k$ sequences $(\gamma_{i,n})_{n \in \N}$, $i \in \{1, \ldots, k\}$, of distinct elements of $\Gamma$ such that for each $\bm{\omega} \in \bracket{k}$ and each permutation $s \in \mathfrak{S}_k$,
    \begin{equation*}
        \lim_{n_{s(1)} \to \infty} \lim_{n_{s(2)} \to \infty} \ldots \lim_{n_{s(k)} \to \infty} T^{\sum_{i=1}^k \omega_i \gamma_{i,n_i}} x_{\bm{0}} = x_{\bm{\omega}}
    \end{equation*}
\end{definition}

Note that Erd\H{o}s progressions correspond to ``degenerate'' Erd\H{o}s cubes where the point $x_{\bm{\omega}}$ depends only on $|\bm{\omega}| = |\{i \in \{1, \ldots, k\} : \omega_i = 1\}|$.

\begin{theorem}[Dynamical reformulation of \Cref{thm-main}] \label{thm: density finite sums dynamical}
    Let $\Gamma=(\Gamma,+)$ be a countably infinite abelian group such that $[\Gamma:\ell\Gamma]<\infty$ for all $\ell\leq k$. Let $(X,T)$ be a topological $\Gamma$-system and $\mu$ an ergodic $T$-invariant Borel probability measure on $X$.
    Suppose $\Phi$ a F{\o}lner sequence in $\Gamma$, $a \in \gen(\mu,\Phi)$, and $E \subseteq X$ is a clopen set with $\mu(E) > 0$.
    Then there exist $x_1, \ldots, x_k \in E$ and $t \in \Gamma$ such that $(T^ta, x_1, \ldots, x_k)$ is an Erd\H{o}s progression.
\end{theorem}

\begin{proof}[Proof that \Cref{thm: density finite sums dynamical} implies \Cref{thm-main}]
    This implication is proved for $\Gamma = \Z$ in \cite[\S 2]{sumsets} for general $\Gamma$ in the case $k = 3$ in \cite[\S 7]{AJ}.
    The general case can be handled by the same argument.
\end{proof}

\begin{theorem}[Dynamical reformulation of \Cref{thm-main1}] \label{thm: sum of k sets dynamical}
    Let $\Gamma=(\Gamma,+)$ be a countably infinite abelian group. Let $(X,T)$ be a topological $\Gamma$-system and $\mu$ an ergodic $T$-invariant Borel probability measure on $X$.
    Suppose $\Phi$ a F{\o}lner sequence in $\Gamma$, $a \in \gen(\mu,\Phi)$, and $E \subseteq X$ is a clopen set with $\mu(E) > 0$.
    Then for every $k \in \N$, there exists a $k$-dimensional Erd\H{o}s cube $(x_{\bm{\omega}})_{\bm{\omega} \in \bracket{k}} \in X^{\bracket{k}}$ such that $x_{\bm{0}} = a$ and $x_{\bm{1}} \in E$.
\end{theorem}

\begin{proof}[Proof that \Cref{thm: sum of k sets dynamical} implies \Cref{thm-main1}]
    The proof for $\Gamma = \Z$ appears in \cite[\S 4.2]{BFRRinfinitesums}, and the same proof works for general $\Gamma$.
\end{proof}

\subsection{Progressive measures}

The existence of Erd\H{o}s progressions in the setting above can be derived from the existence of progressive measures.

\begin{definition} \label{def: progressive}
    Let $(X, T )$ be a topological $\Gamma$-system and $k \in \N$.
    A Borel probability measure $\tau \in \mathcal{M}(X^{k+1})$ is \emph{progressive} if for all open sets $U_1, \ldots, U_k \subseteq X$ with $\tau(X \times U_1 \times \ldots \times U_k) > 0$, there exist infinitely many $\gamma \in \Gamma$ such that
    \begin{equation*}
        \tau \left( \left( X \times U_1 \times \ldots \times U_k \right) \cap \left( T^{-\gamma}U_1 \times \ldots \times T^{-\gamma}U_k \times X\right) \right) > 0.
    \end{equation*}
\end{definition}

\begin{theorem}[Progressive measure reformulation of \Cref{thm: density finite sums dynamical}] \label{thm: progressive measures formulation}
    Let $\Gamma=(\Gamma,+)$ be a countably infinite abelian group such that $[\Gamma:\ell\Gamma]<\infty$ for all $\ell\leq k$.
    Let $(X,T)$ be a topological $\Gamma$-system and $\mu$ an ergodic $T$-invariant Borel probability measure on $X$.
    Suppose $\Phi=(\Phi_N)_{N\in \mathbb{N}}$ a F{\o}lner sequence in $\Gamma$, $a \in \gen(\mu,\Phi)$, and $E \subseteq X$ is a clopen set with $\mu(E) > 0$.
    Then there exists a progressive measures $\tau_a = \delta_a \times \sigma_a \in \M(X^{k+1})$ such that
    \begin{equation*}
        \sigma_{a} \left( T^{-t}E \times \ldots \times T^{-t}E \right) > 0
    \end{equation*}
    for some $t \in \Gamma$.
\end{theorem}

\begin{proof}[Proof that \Cref{thm: progressive measures formulation} implies \Cref{thm: density finite sums dynamical}]
    The case $\Gamma = \Z$ is handled in \cite[\S 3]{sumsets}.
    We explain the general case for completeness.
    Suppose \Cref{thm: progressive measures formulation} holds, and let $\Gamma$ be a countably infinite abelian group such that $[\Gamma : \ell\Gamma] < \infty$ for $\ell \le k$, let $(X, T)$ be a topological $\Gamma$-system, let $\mu$ be a $T$-invariant Borel probability measure on $X$, let $\Phi$ be a F{\o}lner sequence in $\Gamma$, $a \in \gen(\mu,\Phi)$, and $E \subseteq X$ a clopen set with $\mu(E) > 0$.
    By \Cref{thm: progressive measures formulation}, there exists a progressive measure $\tau_a = \delta_a \times \sigma_a$ and $t \in \Gamma$ such that
    \begin{equation*}
        \sigma_a \left( T^{-t}E \times \ldots \times T^{-t}E \right) > 0.
    \end{equation*}
    Combining \cite[Lemma 9.2 and Proposition 9.3]{AJ}, there exists an Erd\H{o}s progression $(a, y_1, \ldots, y_k)$ with $y_1, \ldots, t_k \in T^{-t}E$.
    Hence, putting $x_i = T^ty_i$, we have that $x_1, \ldots, x_k \in E$ and $(T^ta, x_1, \ldots, x_k)$ is an Erd\H{o}s progression.
    This proves \Cref{thm: density finite sums dynamical}.
\end{proof}

\subsection{Erd\H{o}s Cube measures}

Erd\H{o}s cubes can similarly be encoded with measures.

\begin{definition}
    Let $(X, T)$ be a topological $\Gamma$-system and $k \in \N$.
    A Borel probability measure $\rho \in \M \left( X^{\bracket{k}} \right)$ is an \emph{Erd\H{o}s cube measure} if for all open sets $(U_{\bm{\omega}})_{\bm{\omega} \in \bracket{k} \setminus \{\bm{0}\}}$ with $\rho(X \times \prod_{\bm{\omega} \ne \bm{0}} U_{\omega}) > 0$, there are infinite sets $B_1, \ldots, B_k \subseteq \Gamma$ such that
    \begin{equation} \label{eq: cube property}
        \rho \left( \left( \prod_{\bm{\omega} \in \bracket{k}} U_{\bm{\omega}} \right) \cap \bigcap_{i=1}^k \prod_{\bm{\omega} \in \bracket{k}} T^{-b_i} U_{\bm{\omega}+e_i} \right) > 0
    \end{equation}
    for every $(b_1, \ldots, b_k) \in B_1 \times \ldots \times B_k$, where $U_{\bm{0}} = X$, $e_1, \ldots, e_k$ are the standard unit vectors in $\Z^k$, and $U_{\bm{v}} = X$ for $\bm{v} \in \Z^k \setminus \bracket{k}$.
\end{definition}

For example, if $k = 2$, then an Erd\H{o}s cube measure $\rho \in \M \left( X^{\bracket{2}} \right)$ is a measure with the property that if $U_{10}, U_{01}, U_{11} \subseteq X$ are open and $\rho(X \times U_{10} \times U_{01} \times U_{11}) > 0$, then there are infinite sets $B_1, B_2 \subseteq \Gamma$ such that
\begin{multline*}
    \rho \left( (X \times U_{10} \times U_{01} \times U_{11}) \cap T_{\Delta}^{-b_1}(U_{10} \times X \times U_{11} \times X) \right. \\
    \left. \cap T_{\Delta}^{-b_2}(U_{01} \times U_{11} \times X \times X) \right) > 0
\end{multline*}
for every $(b_1, b_2) \in B_1 \times B_2$, where $T_{\Delta} = T \times \ldots \times T$.

For $k = 3$, the expression \eqref{eq: cube property} takes the form
\begin{multline*}
    \rho \left( (X \times U_{100} \times U_{010} \times U_{110} \times U_{001} \times U_{101} \times U_{011} \times U_{111}) \right. \\
    \cap T_{\Delta}^{-b_1}(U_{100} \times X \times U_{110} \times X \times U_{101} \times X \times U_{111} \times X) \\
    \cap T_{\Delta}^{-b_2}(U_{010} \times U_{110} \times X \times X \times U_{011} \times U_{111} \times X \times X) \\
    \left. \cap T_{\Delta}^{-b_3}(U_{001} \times U_{101} \times U_{011} \times U_{111} \times X \times X \times X \times X) \right) > 0.
\end{multline*}

\begin{lemma} \label{lem: cube measures detect cubes}
    Let $(X, T)$ be a topological $\Gamma$-system, $k \in \N$, and $a \in X$.
    Suppose $\rho \in \M \left( X^{\bracket{k}} \right)$ is an Erd\H{o}s cube measure with $\supp(\rho) \subseteq \{a\} \times X^{\bracket{k} \setminus \{\bm{0}\}}$.
    If $(U_{\bm{\omega}})_{\bm{\omega} \in \bracket{k} \setminus \{\bm{0}\}}$ are open subsets of $X$ and $\rho \left( X \times \prod_{\bm{\omega} \ne \bm{0}} U_{\bm{\omega}} \right) > 0$, then there exists a $k$-dimensional Erd\H{o}s cube $(x_{\bm{\omega}})_{\bm{\omega} \in \bracket{k}}$ such that $x_{\bm{0}} = a$ and $x_{\bm{\omega}} \in U_{\bm{\omega}}$ for $\bm{\omega} \in \bracket{k} \setminus \{\bm{0}\}$.
    In particular, if $U \subseteq X$ is open and $\rho \left( X^{\bracket{k} \setminus \{\bm{1}\}} \times U \right) > 0$, then there exists a $k$-dimensional Erd\H{o}s cube $(x_{\bm{\omega}})_{\bm{\omega} \in \bracket{k}}$ such that $x_{\bm{0}} = a$ and $x_{\bm{1}} \in U$.
\end{lemma}

\begin{proof}
    Our proof is inspired by \cite[Proposition 3.2]{sumsets}.
    
    Let $V_{\bm{\omega}}(0) = U_{\bm{\omega}}$ for each $\bm{\omega} \in \bracket{k} \setminus \{\bm{0}\}$.
    We will construct sequences of open subsets $V_{\bm{\omega}}(n) \subseteq X$ and sequences of distinct elements $b_1(n), \ldots, b_k(n) \in \Gamma$ such that for every $n \in \N$,
    \begin{enumerate}
        \item for all $\bm{\omega} \in \bracket{k} \setminus \{\bm{0}\}$, $\overline{V_{\bm{\omega}}(n)} \subseteq V_{\bm{\omega}}(n-1)$ and $\diam \left( V_{\bm{\omega}}(n) \right) \leq \frac{1}{2} \diam \left( V_{\bm{\omega}}(n-1) \right)$;
        \item $\rho \left( \{a\} \times \prod_{\bm{\omega} \in \bracket{k} \setminus \{\bm{0}\}} V_{\bm{\omega}}(n) \right) > 0$;
        \item for all $i \in \{1, \ldots, k\}$,
        \begin{equation*}
            T^{b_i(n)} a \in V_{e_i}(n-1);
        \end{equation*}
        \item for all $\bm{\omega} \in \bracket{k} \setminus \{\bm{0}\}$ and all $i \in \{1, \ldots, k\}$,
        \begin{equation*}
            V_{\bm{\omega}}(n) \subseteq T^{-b_i(n)} V_{\bm{\omega} + e_i}(n-1);
        \end{equation*}
        \item for all $i \in \{1, \ldots, k\}$, $b_i(n) \notin \{b_i(j) : j < n\}$.
    \end{enumerate}

    First, the hypothesis of the lemma gives (2) for $n=0$.
    Since $\rho$ is an Erd\H{o}s cube measure, it follows that there exist $b_1(1), \ldots, b_k(1) \in \Gamma$ such that
    \begin{equation*}
        \rho \left( \left( \{a\} \times \prod_{\bm{\omega} \ne \bm{0}} V_{\bm{\omega}}(0) \right) \cap \bigcap_{i=1}^k T_{\Delta}^{-b_i(1)} \prod_{\bm{\omega} \in \bracket{k}} V_{\bm{\omega}+e_i}(0) \right) > 0.
    \end{equation*}
    By regularity of the measure $\rho$, we may then find open sets $V_{\bm{\omega}}(1) \subseteq X$ with $\diam(V_{\bm{\omega}}(1)) \leq \frac{1}{2} \diam(V_{\bm{\omega}}(0))$ such that
    \begin{equation} \label{eq: inclusion in intersection}
        \{a\} \times \prod_{\bm{\omega} \ne \bm{0}} \overline{V_{\bm{\omega}}(1)} \subseteq \left( \{a\} \times \prod_{\bm{\omega} \ne \bm{0}} V_{\bm{\omega}}(0) \right) \cap \bigcap_{i=1}^k T_{\Delta}^{-b_i(1)} \prod_{\bm{\omega} \in \bracket{k}} V_{\bm{\omega}+e_i}(0)
    \end{equation}
    and
    \begin{equation*}
        \rho \left( \{a\} \times \prod_{\bm{\omega} \ne \bm{0}} V_{\bm{\omega}}(1) \right) > 0.
    \end{equation*}
    This final inequality gives (2) for $n=1$.
    Property (1) holds for $n=1$ by looking at the first part of the inclusion in \eqref{eq: inclusion in intersection} coordinate-wise and (3) and (4) also hold for $n=1$ by looking at the second part of the inclusion coordinate-wise.
    Property (5) is vacuous for $n=1$, so we have completed the base case.

    The induction step is exactly the same.
    The only additional thing to check is property (5), but this can be achieved since the assumption of being an Erd\H{o}s cube measure means that \eqref{eq: cube property} holds for all $b_1, \ldots, b_k$ coming from \emph{infinite} sets $B_1, \ldots, B_k \subseteq \Gamma$, so in particular $B_i \setminus \{b_i(j) : j < n\} \ne \emptyset$.

    \bigskip

    Let us now use the constructed sequences to conclude.
    By (1), the intersection $\bigcap_{n \in \N}V_{\bm{\omega}}(n)$ contains a single point $x_{\bm{\omega}}$ for $\bm{\omega} \in \bracket{k} \setminus \{\bm{0}\}$.
    Moreover, $x_{\bm{\omega}} \in V_{\bm{\omega}}(0) = U_{\bm{\omega}}$.
    Let $x_{\bm{0}} = a$.
    We claim that $(x_{\bm{\omega}})_{\bm{\omega} \in \bracket{k}}$ is a $k$-dimensional Erd\H{o}s cube.
    Fix $\bm{\omega} \in \bracket{k} \setminus \{\bm{0}\}$.
    Suppose $i_1, \ldots, i_\ell$ is a permutation of $\{i \in \{1, \ldots, k\} : \omega_i = 1\}$ and $n_{i_1} < \ldots < n_{i_\ell}$.
    Then
    \begin{align*}
        T^{\sum_{j=1}^\ell b_{i_j}(n_{i_j})} a & = T^{b_{i_1}(n_{i_1})} \cdot \ldots \cdot T^{b_{i_\ell}(n_{i_\ell})} a \\
        & \in T^{b_{i_1}(n_{i_1})} \cdot \ldots \cdot T^{b_{i_{\ell-1}}(n_{i_{\ell-1}})} V_{e_{i_\ell}}(n_\ell-1) \\
        & \subseteq T^{b_{i_1}(n_{i_1})} \cdot \ldots \cdot T^{b_{i_{\ell-1}}(n_{i_{\ell-1}})} V_{e_{i_\ell}}(n_{\ell-1}) \\
        & \subseteq T^{b_{i_1}(n_{i_1})} \cdot \ldots \cdot T^{b_{i_{\ell-2}}(n_{i_{\ell-2}})} V_{e_{i_\ell} + e_{i_{\ell-1}}}(n_{\ell-1}-1) \\
        & \subseteq T^{b_{i_1}(n_{i_1})} \cdot \ldots \cdot T^{b_{i_{\ell-2}}(n_{i_{\ell-2}})} V_{e_{i_\ell} + e_{i_{\ell-1}}}(n_{\ell-2}) \\
        & \qquad \vdots \\
        & \subseteq V_{\bm{\omega}}(n_{i_1}-1),
    \end{align*}
    since $\bm{\omega} = \sum_{j=1}^\ell e_{i_j}$.
    Hence,
    \begin{equation*}
        \lim_{n_{i_1} \to \infty} \ldots \lim_{n_{i_\ell} \to \infty} T^{\sum_{j=1}^\ell b_{i_j}(n_{i_j})} a = x_{\bm{\omega}}.
    \end{equation*}
    This proves that $(x_{\bm{\omega}})_{\bm{\omega} \in \bracket{k}}$ is an Erd\H{o}s cube (see \Cref{def: erdos cube}).
\end{proof}

\begin{theorem}[Erd\H{o}s cube measure reformulation of \Cref{thm: sum of k sets dynamical}] \label{thm: Erdos cube measures formulation}
    Let $\Gamma=(\Gamma,+)$ be a countably infinite abelian group and $k \in \N$.
    Let $(X,T)$ be a topological $\Gamma$-system and $\mu$ an ergodic $T$-invariant Borel probability measure on $X$.
    Suppose $\Phi=(\Phi_N)_{N\in \mathbb{N}}$ a F{\o}lner sequence in $\Gamma$, $a \in \gen(\mu,\Phi)$, and $E \subseteq X$ is a clopen set with $\mu(E) > 0$.
    Then there exists an Erd\H{o}s cube measure $\rho_a \in \M(X^{\bracket{k}})$ such that $\supp(\rho_a) \subseteq \{a\} \times X^{\bracket{k} \setminus \{\bm{0}\}}$ and
    \begin{equation*}
        \rho_a \left( X^{\bracket{k} \setminus \{\bm{1}\}} \times E \right) > 0.
    \end{equation*}
\end{theorem}

\begin{proof}[Proof that \Cref{thm: Erdos cube measures formulation} implies \Cref{thm: sum of k sets dynamical}]
    This follows immediately from \Cref{lem: cube measures detect cubes}.
\end{proof}

\section{Proof of \Cref{thm-main}} \label{sec:proof of density FS}

In this section, we prove \Cref{thm-main}.
As shown in \Cref{sec:dynamical reformulation}, it suffices to prove \Cref{thm: progressive measures formulation}.
Our main task, then, is to give a construction of suitable progressive measures.

Throughout this section, we fix a countably infinite abelian group $\Gamma=(\Gamma,+)$ satisfying $[\Gamma:\ell\Gamma]<\infty$ for all $\ell\leq k$, equivalently $[\Gamma:k!\Gamma]<\infty$, topological dynamical $\Gamma$-system $(X,T)$, an ergodic $T$-invariant Borel probability measure $\mu \in \M(X)$, a F{\o}lner sequence $\Phi=(\Phi_N)_{N\in \mathbb{N}}$ in $\Gamma$, a point $a \in \gen(\mu,\Phi)$, and a clopen set $E \subseteq X$ with $\mu(E) > 0$, as in the hypothesis of \Cref{thm: progressive measures formulation}.
Our goal is to construct a progressive measure satisfying the conclusion of \Cref{thm: progressive measures formulation}.

By \cite[Theorem 5.11]{candela-szegedy-cubic} and \cite[Proposition 5.4]{CGSSextensions}, let $(Z, R)$ be a minimal, distal, and uniquely ergodic system (with unique invariant measure $\mu_Z)$ of order $k-1$ such that the $(k-1)$-step Host--Kra factor of $(X, \mu, T)$ is measurably isomorphic to $(Z, \mu_Z, R)$.
By \cite[Proposition 3.7]{cm}, we can replace $(X,\mu,T)$ by an extension so that the factor map $\pi \colon X \to Z$ is continuous.
(Strictly speaking, \cite[Proposition 3.7]{cm} handles the case that $Z$ is the Kronecker factor, corresponding to $k=2$.
However, the method applies with any topologically distal factor; see the remark in \cite[Lemma 3.5]{cm}.)
We will first construct a progressive measure for the factor $(Z, R)$ and then lift to $(X,T)$.

\subsection{Defining a measure}

Let $c = \pi(a)$.
The system $\left( Z^k, R \times R^2 \times \ldots \times R^k \right)$ is a system of order $k$.
Therefore, by \Cref{thm:minimal finite order systems are u.e.}, there is a unique invariant measure on the orbit-closure
\begin{equation*}
    \mathcal{O}(c) = \overline{\{(R^{\gamma}c, R^{2\gamma}c, \ldots, R^{k\gamma}c) : \gamma \in \Gamma\}}.
\end{equation*}
Denote this measure by $\xi$.

\begin{lemma} \label{lem:marginals}
    For each $\ell \in [k]$, let $p_{\ell} \colon Z^k \to Z$ be the projection onto the $\ell$th coordinate.
    Then $(p_\ell)_*\xi \leq [\Gamma:\ell\Gamma] \cdot \mu_Z$.
\end{lemma}
\begin{proof}
    Let $m = (p_\ell)_*\xi \in \M(Z)$.
    By construction, $m$ is $R^\ell$-invariant.
    Consider the measure $m' \in \M(Z)$ defined by
    \begin{equation*}
        m'(B) = \frac{1}{[\Gamma:\ell\Gamma]} \sum_{t \in \Gamma/\ell\Gamma} m(R^{-t}B)
    \end{equation*}
    for Borel sets $B \subseteq Z$.
    (Note that, because $m$ is $R^{\ell}$-invariant, the quantity $m(R^{-t}A)$ does not depend on the representative of the coset $t$, so $m'$ is well-defined.)
    Then $m'$ is $R$-invariant, so by unique ergodicity of $(Z, R)$, we conclude that $m'=\mu_Z$.
    Hence,
    \begin{equation*}
        m(B) = [\Gamma:\ell\Gamma] \cdot \mu_Z(B) - \sum_{t \in (\Gamma/\ell\Gamma) \setminus \{0\}} m(R^{-t}B) \leq [\Gamma:\ell\Gamma] \cdot \mu_Z(B).
    \end{equation*}
\end{proof}

Define the measure $\sigma_a \in \M(X^k)$ by
\begin{equation*}
    \int_{X^k} f_1 \otimes \ldots \otimes f_k~d\sigma_a = \int_{Z^k} \E[f_1 \mid Z] \otimes \ldots \otimes \E[f_1 \mid Z]~d\xi
\end{equation*}
for $f_1, \ldots, f_k \in C(X)$.
By \Cref{lem:marginals}, the measure $\sigma_a$ is well-defined, i.e., does not depend on the choice of representatives of the conditional expectations $\E[f_i \mid Z]$.
Let $\tau_a = \delta_a \times \sigma_a \in \M(X^{k+1})$.

\subsection{Verifying part of the conclusion of \Cref{thm: progressive measures formulation}}

We will show in \Cref{sec:progressive} that $\tau_a$ is a progressive measure.
Before doing so, we verify in this subsection that $\sigma_a$ satisfies the desired property
\begin{equation*}
    \sigma_a(T^{-t}E \times \ldots \times T^{-t}E) > 0
\end{equation*}
for some $t \in \Gamma$ (see \Cref{cor:positive mass to shift of E} below).
As in \cite{sumsets}, this is done by comparing with the Furstenberg joining.

\begin{definition}[Furstenberg joining]
    The \emph{Furstenberg joining} is the measure $\lambda \in \M(X^{k+1})$ defined by
    \begin{equation*}
        \int_{X^{k+1}} f_0 \otimes f_1 \otimes \ldots \otimes f_k~d\lambda = \lim_{N \to \infty} \frac{1}{|\Theta_N|} \sum_{\gamma \in \Theta_N} \int_X f_0 \cdot T^{\gamma}f_1 \cdot \ldots \cdot T^{k\gamma} f_k~d\mu
    \end{equation*}
    for $f_0, f_1, \ldots, f_k \in C(X)$, where $\Theta$ is any F{\o}lner sequence in $\Gamma$.
\end{definition}

The Furstenberg joining is well-defined (i.e., the limit exists and does not depend on the choice of F{\o}lner sequence) by \cite{austin2016nonconventional, zorin2016norm}.

The Host--Kra seminorms for arbitrary countable abelian $\Gamma$-systems were defined in \cite[Appendix A]{bergelson2010inverse}.

\begin{lemma}[Seminorm control]\label{seminorm-control}
Let $k\geq 1$. Let $\Gamma$ be a countable abelian group and suppose that $[\Gamma:k!\Gamma]<\infty$. Let $(X,\mu,T)$ be a measure-preserving $\Gamma$-system and $f_1,\dots,f_k$ bounded functions on $X$ with $\|f_j\|_\infty \leq 1$ for all $j=1,\dots,k$. Then, for any F\o lner sequence $\Theta$ in $\Gamma$, we have
$$\limsup_{N\rightarrow\infty} \left\|\frac{1}{|\Theta_N|}\sum_{\gamma\in \Phi_N} \prod_{i=1}^k T^{i\gamma} f_i\right\|_2 \leq [\Gamma:k!\Gamma]\cdot \min_{1\leq i \leq k} \|f_i\|_{U^{k}(\XX)}$$
\end{lemma}

\begin{proof}
    The proof is the same as in \cite[Theorem 12.1]{host2005nonconventional} where in that proof the coefficient $\ell$ is the index $[\Z:\ell\Z]$. The proof is an iterated application of the van der Corput inequality which was generalized to arbitrary countable abelian groups in \cite[Proposition 6.4.27]{bergelsondiophantine}.
\end{proof}
\begin{theorem} \label{thm:average over shifts Furstenberg joining}
    For $f_1, \ldots, f_k \in C(X)$ and any F{\o}lner sequence $\Theta$ in $\Gamma$,
    \begin{equation} \label{eq:average over translates = furstenberg}
        \lim_{N \to \infty} \frac{1}{|\Theta_N|} \sum_{\gamma \in \Phi_N} \int_{X^{k+1}} T_{\Delta}^\gamma \left( 1 \otimes f_1 \otimes \ldots \otimes f_k\right)~d\tau_a = \int_{X^{k+1}} 1 \otimes f_1 \otimes \ldots \otimes f_k~d\lambda.
    \end{equation}
\end{theorem}

\begin{proof}
    By \Cref{seminorm-control}, we have
    \begin{equation*}
        \int_{X^{k+1}} 1 \otimes f_1 \otimes \ldots \otimes f_k~d\lambda = \lim_{N \to \infty} \frac{1}{|\Theta_N|} \sum_{\gamma \in \Theta_N} \int_Z R^{\gamma} \varphi_1 \cdot \ldots \cdot R^{k\gamma} \varphi_k~d\mu_Z,
    \end{equation*}
    where $\varphi_i = \E[f_i \mid Z]$.
    Similarly, the left-hand side of \eqref{eq:average over translates = furstenberg} can be rewritten as
    \begin{equation*}
        \lim_{N \to \infty} \frac{1}{|\Theta_N|} \sum_{\gamma \in \Phi_N} \int_{Z^{k+1}} R^\gamma \varphi_1 \otimes \ldots \otimes R^{\gamma} \varphi_k~d\xi.
    \end{equation*}
    Therefore, it suffices to prove: for all $\psi_1, \ldots, \psi_k \in L^{\infty}(\mu_Z)$,
    \begin{equation} \begin{split} \label{eq:interchange of limits}
        \lim_{N \to \infty} & \frac{1}{|\Theta_N|} \sum_{\gamma \in \Theta_N} \int_{Z^{k+1}} R^\gamma \psi_1 \otimes \ldots \otimes R^{\gamma} \psi_k~d\xi \\
         & = \lim_{N \to \infty} \frac{1}{|\Theta_N|} \sum_{\gamma \in \Theta_N} \int_Z R^{\gamma} \psi_1 \cdot \ldots \cdot R^{k\gamma} \psi_k~d\mu_Z.
    \end{split} \end{equation}
    By a standard approximation argument, it suffices to consider the case $\psi_1, \ldots, \psi_k \in C(Z)$ in \eqref{eq:interchange of limits}.
    Consider the $(\Gamma \times \Gamma)$-action
    \begin{equation*}
        S^{\bm{\gamma}} = R^{\gamma_1+\gamma_2} \times R^{\gamma_1+2\gamma_2} \times \ldots \times R^{\gamma_1+k\gamma_2}.
    \end{equation*}
    Let
    \begin{equation*}
        Y = \overline{\left\{ S^{\bm{\gamma}}(c,\ldots,c) : \bm{\gamma} \in \Gamma \times \Gamma\right\}} \subseteq Z^k.
    \end{equation*}
    As a $(\Gamma \times \Gamma)$-system, $(Y, S)$ is minimal of order $k$ and hence uniquely ergodic by \Cref{thm:minimal finite order systems are u.e.}.
    Let $\nu$ be the unique $S$-invariant measure on $Y$.
    Unfolding the left-hand side of \eqref{eq:interchange of limits},
    \begin{align*}
         \lim_{N \to \infty} & \frac{1}{|\Theta_N|} \sum_{\gamma \in \Theta_N} \int_{Z^{k+1}} R^\gamma \psi_1 \otimes \ldots \otimes R^{\gamma} \psi_k~d\xi \\
          & = \lim_{N \to \infty} \frac{1}{|\Theta_N|} \sum_{\gamma_1 \in \Theta_N}  \lim_{M \to \infty} \frac{1}{|\Theta_M|} \sum_{\gamma_2 \in \Theta_M} \psi_1(R^{\gamma_1+\gamma_2}c) \cdot \psi_2(R^{\gamma_1 + 2\gamma_2}c) \cdot \ldots \cdot \psi_k(R^{\gamma_1+k\gamma_2}c) \\
          & = \int_{Z^k} \psi_1 \otimes \ldots \otimes \psi_k~d\nu.
    \end{align*}
    Similarly, using unique ergodicity of $(Z, R)$ to unfold the right-hand side of \eqref{eq:interchange of limits}, we have
    \begin{align*}
        \lim_{N \to \infty} & \frac{1}{|\Theta_N|} \sum_{\gamma \in \Theta_N} \int_Z R^{\gamma} \psi_1 \cdot \ldots \cdot R^{k\gamma} \psi_k~d\mu_Z \\
        & = \lim_{N \to \infty} \frac{1}{|\Theta_N|} \sum_{\gamma_2 \in \Theta_N} \lim_{M \to \infty} \frac{1}{|\Theta_M|} \sum_{\gamma_1 \in \Theta_M} \psi_1(R^{\gamma_1+\gamma_2}c) \cdot \psi_2(R^{\gamma_1+2\gamma_2}c) \cdot \ldots \cdot \psi_k(R^{\gamma_1+k\gamma_2}c) \\
        & = \int_{Z^k} \psi_1 \otimes \ldots \otimes \psi_k~d\nu.
    \end{align*}
\end{proof}

\begin{corollary} \label{cor:positive mass to shift of E}
    There exists $t \in \Gamma$ such that
    \begin{equation*}
        \sigma_a \left( T^{-t}E \times \ldots \times T^{-t}E \right) > 0.
    \end{equation*}
\end{corollary}

\begin{proof}
    Let $\Theta$ be a F{\o}lner sequence in $\Gamma$.
    By \Cref{thm:average over shifts Furstenberg joining} and Szemer\'{e}di's theorem (see, e.g., \cite[Theorem B]{austin2016nonconventional}),
    \begin{align*}
        \lim_{N \to \infty} \frac{1}{|\Theta_N|} & \sum_{t \in \Theta_N} \sigma_a\left( T^{-t}E \times \ldots \times T^{-t}E \right) \\
        & = \lim_{N \to \infty} \frac{1}{|\Theta_N|} \sum_{\gamma \in \Theta_N} \int_{X^{k+1}} T_{\Delta}^t \left( 1 \otimes 1_E \otimes \ldots \otimes 1_E \right)~d\tau_a \\
        & = \int_{X^{k+1}} 1 \otimes 1_E \otimes \ldots \otimes 1_E~d\lambda \\
        & = \lim_{N \to \infty} \frac{1}{|\Theta_N|} \sum_{\gamma \in \Theta_N} \mu \left( T^{-\gamma}E \cap \ldots \cap T^{-k\gamma}E \right) > 0.
    \end{align*}
\end{proof}

\subsection{Proving progressiveness} \label{sec:progressive}

It remains to prove that $\tau_a$ is progressive (\Cref{def: progressive}).
We will prove progressiveness of $\tau_a$ using a functional version: if $f_1, \ldots, f_k \colon X \to [0,1]$ are continuous functions such that
\begin{equation*}
    \int_{X^{k+1}} 1 \otimes f_1 \otimes \ldots \otimes f_k~d\tau_a > 0,
\end{equation*}
then there are infinitely many $\gamma \in \Gamma$ such that
\begin{equation*}
    \int_{X^{k+1}} T^{\gamma}f_1 \otimes (f_1 T^{\gamma}f_2) \otimes \ldots \otimes (f_{k-1} T^{\gamma}f_k) \otimes f_k~d\tau_a > 0.
\end{equation*}
We will prove this by averaging over $\gamma \in \Gamma$.

The first step for handling the average over $\gamma \in \Gamma$ is to transfer functions between coordinates.

\begin{lemma} \label{lem:coordinate swap}
    Let $f_1, \ldots, f_k \in C(X)$.
    Then for any F{\o}lner sequence $\Theta$ in $\Gamma$,
    \begin{equation*}
        \lim_{N \to \infty} \left\| \frac{1}{|\Theta_N|} \sum_{\gamma \in \Theta_N} T^{\gamma}f_1 \otimes \ldots \otimes T^{\gamma}f_k \otimes 1 - \frac{1}{|\Theta_N|} \sum_{\gamma \in \Theta_N}1 \otimes T^{\gamma}f_1 \otimes \ldots \otimes T^{\gamma}f_k \right\|_{L^2(\tau_a)} = 0.
    \end{equation*}
\end{lemma}

\begin{proof}
    This is established for $\Gamma = \Z$ in \cite{sumsets}.
    An inspection of the proof shows that the arguments directly extend to our setting.
\end{proof}

The second step for handling the average is to prove a recurrence result.

\begin{lemma} \label{lem:diagonal recurrence}
    Let $\sigma$ be a $(T \times T^2 \times \ldots \times T^k)$-invariant measure on $X^k$, and let $A_1, \ldots, A_k \subseteq X$ be measurable sets.
    Put $A = A_1 \times \ldots \times A_k$.
    Then for any F{\o}lner sequence $\Theta$ in $\Gamma$,
    \begin{equation*}
        \liminf_{N \to \infty} \frac{1}{|\Theta_N|} \sum_{\gamma \in \Theta_N} \sigma \left( A \cap T_{\Delta}^{\gamma}A \right) > 0.
    \end{equation*}
\end{lemma}

\begin{proof}
    A version of this statement (for $\Z$-actions) appears in \cite[Theorem 6.3]{sumsets}.
    We follow a similar approach, using \Cref{lem:Folner lifting} to overcome some technical hurdles involving multiplication by $k!$.
    
    The basic idea is to bound the correlation sequence $\sigma(A \cap T_{\Delta}^{-\gamma}A)$ from below by a multicorrelation sequence for the action $S = T \times T^2 \times \ldots \times T^k$ for which $\sigma$ is invariant in order to apply Szemer\'{e}di's theorem.
    
    Observe that
    \begin{align*}
        T_{\Delta}^{-k!\gamma}A & = \bigcap_{i=1}^k X^{i-1} \times T^{-k!\gamma} A_i \times X^{k-i} \\
        & = \bigcap_{i=1}^k S^{-\frac{k!}{i}\gamma} \left( X^{i-1} \times A_i \times X^{k-i} \right) \\
        & \supseteq \bigcap_{i=1}^k S^{- \frac{k!}{i}\gamma} A.
    \end{align*}

    Now let $q \colon \Gamma \to k! \Gamma$ be the map $q(\gamma) = k!\gamma$.
    Let $\Theta'_N = \Theta_N \cap k!\Gamma$ to obtain a F{\o}lner sequence in $k!\Gamma$.
    By \Cref{lem:Folner lifting}, there exists a F{\o}lner sequence $\tilde{\Theta}$ in $\Gamma$ such that $q(\tilde{\Theta}_N) = \Theta'_N$ and $q$ is constant-to-one on $\tilde{\Theta}_N$.
    We then have
    \begin{align*}
        \liminf_{N \to \infty} \frac{1}{|\Theta_N|} \sum_{\gamma \in \Theta_N} \sigma(A \cap T_{\Delta}^{\gamma}A) 
        & \geq \frac{1}{[\Gamma : k!\Gamma]} \liminf_{N \to \infty} \frac{1}{|\Theta'_N|} \sum_{\gamma' \in \Theta'_N} \sigma(A \cap T_{\Delta}^{\gamma'}A) \\
        & = \frac{1}{[\Gamma : k!\Gamma]} \liminf_{N \to \infty} \frac{1}{|\tilde{\Theta}_N|} \sum_{\gamma \in \tilde{\Theta}_N} \sigma(A \cap T_{\Delta}^{k!\gamma} A) \\
        & \geq \frac{1}{[\Gamma : k!\Gamma]} \liminf_{N \to \infty} \frac{1}{|\tilde{\Theta}_N|} \sum_{\gamma \in \tilde{\Theta}_N} \sigma \left( A \cap \bigcap_{i=1}^k S^{-\frac{k!}{i}\gamma}A \right),
    \end{align*}
    and this final quantity is positive by Szemer\'{e}di's theorem (see, e.g., \cite[Theorem B]{austin2016nonconventional}).
\end{proof}

\begin{theorem}
    $\tau_a$ is progressive.
\end{theorem}

\begin{proof}
    Recall $\tau_a = \delta_a \times \sigma_a$.
    Let $f_1, \ldots, f_k \colon X \to [0,1]$ be continuous functions with
    \begin{equation*}
        \int_{X^k} \bigotimes_{i=1}^k f_i~d\sigma_a = \int_{X^{k+1}} 1 \otimes f_1 \otimes \ldots \otimes f_k~d\tau_a > 0.
    \end{equation*}
    Let $\Theta$ be a F{\o}lner sequence in $\Gamma$.
    Applying \Cref{lem:coordinate swap} and then \Cref{lem:diagonal recurrence}, we have
    \begin{align*}
        \liminf_{N \to \infty} \frac{1}{|\Theta_N|} & \sum_{\gamma \in \Theta_N} \int_{X^{k+1}} T^{\gamma}f_1 \otimes (f_1 T^{\gamma}f_2) \otimes \ldots \otimes (f_{k-1} T^{\gamma}f_k) \otimes f_k~d\tau_a \\
         & = \liminf_{N \to \infty} \frac{1}{|\Theta_N|} \sum_{\gamma \in \Theta_N} \int_{X^{k+1}} 1 \otimes (f_1 T^{\gamma}f_1) \otimes \ldots \otimes (f_{k-1} T^{\gamma}f_{k-1}) \otimes (f_k T^{\gamma}f_k)~d\tau_a \\
         & = \liminf_{N \to \infty} \frac{1}{|\Theta_N|} \int_{X^k} F \otimes T_{\Delta}^{\gamma}F~d\sigma_a > 0,
    \end{align*}
    where $F = \bigotimes_{i=1}^k f_i$.
    In particular, there are infinitely many $\gamma \in \Gamma$ such that
    \begin{equation*}
        \int_{X^{k+1}} T^{\gamma}f_1 \otimes (f_1 T^{\gamma}f_2) \otimes \ldots \otimes (f_{k-1} T^{\gamma}f_k) \otimes f_k~d\tau_a > 0.
    \end{equation*}
    This proves that $\tau_a$ is progressive.
\end{proof}

\section{Proof of \Cref{thm-main1}} \label{sec:proof sum of k sets}

As in \Cref{sec:proof of density FS}, we begin by fixing dynamical data.
Throughout this section, $\Gamma=(\Gamma,+)$ is a fixed (but arbitrary) countably infinite abelian group, $(X, T)$ is a topological dynamical $\Gamma$-system with an ergodic measure $\mu$, $\Phi$ is a F{\o}lner sequence in $\Gamma$, $a \in \gen(\mu,\Phi)$, and $E \subseteq X$ is a clopen set with $\mu(E) > 0$.
We assume that $(Z,R)$ is a minimal, distal, and uniquely ergodic system of order $k-1$ such that the $(k-1)$-step Host--Kra factor of $(X, \mu, T)$ is measurably isomorphic to $(Z, \mu_Z, R)$, where $\mu_Z$ is the unique $R$-invariant Borel probability measure on $Z$, and we assume that $\pi \colon X \to Z$ is a continuous factor map.

Our goal is to construct an Erd\H{o}s cube measure satisfying the conclusion of \Cref{thm: Erdos cube measures formulation}.

\subsection{Defining a measure}

Let $c = \pi(a)$.
We consider the face action on $Z^{\bracket{k}}$.
That is, for $\bm{\gamma} = (\gamma_1, \ldots, \gamma_k) \in \Gamma^k$, let
\begin{equation*}
    S^{\bm{\gamma}} (z_{\bm{\omega}})_{\bm{\omega} \in \bracket{k}} = \left( R^{\bm{\omega} \cdot \bm{\gamma}} z_{\bm{\omega}} \right)_{\bm{\omega}\in\bracket{k}}
\end{equation*}
for $(z_{\bm{\omega}})_{\bm{\omega} \in \bracket{k}} \in Z^{\bracket{k}}$.
By \Cref{thm:minimal finite order systems are u.e.}, let $\zeta$ be the unique $S$-invariant measure on
\begin{equation*}
    Q_c^{\bracket{k}} = \overline{\{S^{\bm{\gamma}}(c,c,\ldots,c) : \bm{\gamma} \in \Gamma^k\}} = \{(z_{\bm{\omega}})_{\bm{\omega}\in\bracket{k}} \in Q^{\bracket{k}} : z_{\bm{0}} = c\}.
\end{equation*}

\begin{lemma} \label{lem:cube marginals}
    For each $\bm{\omega} \in \bracket{k}$, let $p_{\bm{\omega}} \colon Z^{\bracket{k}} \to Z$ be the projection onto the $\bm{\omega}$ coordinate.
    Then $(p_{\bm{\omega}})_* \zeta = \mu_Z$ for every $\bm{\omega} \in \bracket{k} \setminus \{\bm{0}\}$.
\end{lemma}

\begin{proof}
    Let $\bm{\omega} \in \bracket{k} \setminus \{\bm{0}\}$, and let $m = (p_{\bm{\omega}})_* \zeta$.
    Since $\bm{\omega} \ne \bm{0}$, let $i \in [k]$ with $\omega_i = 1$.
    Taking elements of the form $\gamma e_i \in \Gamma^k$ for $\gamma \in \Gamma$, we have
    \begin{equation*}
        m(R^{-\gamma}B) = (p_{\bm{\omega}})_* \zeta(S^{-\gamma e_i}B) = (p_{\bm{\omega}})_* \zeta(B) = m(B),
    \end{equation*}
    so $m$ is $R$-invariant.
    But $(Z,R)$ is uniquely ergodic, so $m=\mu_Z$.
\end{proof}

Let $\zeta_* \in \M(Z^{\bracket{k} \setminus \{\bm{0}\}})$ be the measure such that $\zeta = \delta_c \times \zeta_*$.
We then define a measure $\rho_a \in \M(X^{\bracket{k}})$ by
\begin{equation*}
    \int_{X^{\bracket{k}}} \bigotimes_{\bm{\omega} \in \bracket{k}} f_{\bm{\omega}}~d\rho_a = f_{\bm{0}}(a) \int_{Z^{\bracket{k} \setminus \{\bm{0}\}}} \bigotimes_{\bm{\omega} \ne \bm{0}} \E[f_{\bm{\omega}} \mid Z]~d\zeta_*
\end{equation*}
for $f_{\bm{\omega}} \in C(X)$.
By \Cref{lem:cube marginals}, $\rho_a$ is well-defined.

\begin{proposition}
    $\rho_a(X^{\bracket{k}\setminus\{\bm{1}\}} \times E) > 0$.
\end{proposition}

\begin{proof}
    The function $1_E$ is continuous, so by definition
    \begin{equation*}
        \rho_a(X^{\bracket{k}\setminus\{\bm{1}\}} \times E) = \int_{X^{\bracket{k}}} \left( \bigotimes_{\bm{\omega} \in \bracket{k} \setminus \{\bm{1}\}} 1 \right) \otimes 1_E~d\rho_a = \int_{Z^{\bracket{k} \setminus \{\bm{0}\}}} \left( \bigotimes_{\bm{\omega} \ne \bm{0}, \bm{1}} 1 \right) \otimes \E[1_E \mid Z]~d\zeta_*,
    \end{equation*}
    and by \Cref{lem:cube marginals}, this is equal to
    \begin{equation*}
        \int_Z \E[1_E \mid Z]~d\mu_Z = \mu(E) > 0.
    \end{equation*}
\end{proof}

\subsection{Proving the measure is an Erd\H{o}s cube measure}

We now turn to proving that $\rho_a$ is an Erd\H{o}s cube measure.
The strategy is broadly similar to the strategy used in \Cref{sec:progressive} above to establish progressiveness of $\tau_a$.

\begin{lemma} \label{lem:equivalent to face action}
    For each $\bm{\omega} \in \bracket{k} \setminus \{\bm{0}\}$, let $f_{\bm{\omega}} \in C(X)$.
    Let $f_{\bm{0}} = 1$ and $f_{\bm{v}} = 1$ for $\bm{v} \in \Z^k \setminus \bracket{k}$.
    Then for every F{\o}lner sequence $\Theta$ in $\Gamma^k$,
    \begin{equation} \label{eq:cube coordinate swap}
        \lim_{N \to \infty} \left\| \frac{1}{|\Theta_N|} \sum_{\bm{\gamma} \in \theta_N} \prod_{i=1}^k \bigotimes_{\bm{\omega} \in \bracket{k}} T^{\gamma_i} f_{\bm{\omega}+e_i} - \frac{1}{|\Theta_N|} \sum_{\bm{\gamma} \in \theta_N} \prod_{i=1}^k T_i^{\gamma_i} \bigotimes_{\bm{\omega} \in \bracket{k}}f_{i,\bm{\omega}} \right\|_{L^2(\rho_a)} = 0,
    \end{equation}
    where $T_i \colon \Gamma \acts X^{\bracket{k}}$ is the $i$th face action (see \Cref{def:parallelepiped group}) and
    \begin{equation*}
        f_{i,\bm{\omega}} = \begin{cases}
            f_{\bm{\omega}}, & \text{if}~\omega_i = 1; \\
            1, & \text{if}~\omega_i=0.
        \end{cases}
    \end{equation*}
\end{lemma}

\begin{proof}
    Arguing similarly to \Cref{lem:coordinate swap} above (justified in \cite{sumsets}), the averages in \eqref{eq:cube coordinate swap} are controlled by the order $(k-1)$ factor $Z$, so we can reduce to showing the result in this factor.

    Assume that $h_{\bm{\omega}} \in C(Z)$ for $\bm{\omega} \in \bracket{k} \setminus \{\bm{0}\}$, and set $h_{\bm{0}} = 1$ and $h_{\bm{v}} = 1$ for $\bm{v} \in \Z^k \setminus \bracket{k}$.
    By construction,
    \begin{align*}
        \lim_{N \to \infty} & \left\| \frac{1}{|\Theta_N|} \sum_{\bm{\gamma} \in \theta_N} \prod_{i=1}^k \bigotimes_{\bm{\omega} \in \bracket{k}} T^{\gamma_i} h_{\bm{\omega}+e_i} - \frac{1}{|\Theta_N|} \sum_{\bm{\gamma} \in \theta_N} \prod_{i=1}^k T_i^{\gamma_i} \bigotimes_{\bm{\omega} \in \bracket{k}} h_{i,\bm{\omega}} \right\|_{L^2(\zeta)}^2 \\
        & = \lim_{N \to \infty} \lim_{M \to \infty} \frac{1}{|\Theta_M|} \sum_{\bm{\delta} \in \Theta_M} \Bigg| \frac{1}{|\Theta_N|} \sum_{\bm{\gamma} \in \Theta_N} \prod_{i=1}^k \prod_{\bm{\omega} \in \bracket{k}} h_{\bm{\omega}+e_i}(T^{\gamma_i+\bm{\omega}\cdot\bm{\delta}}c) \\
        & \qquad \qquad \qquad \qquad \qquad \qquad - \frac{1}{|\Theta_N|} \sum_{\bm{\gamma} \in \Theta_N} \prod_{i=1}^k \prod_{\bm{\omega} \in \bracket{k}} h_{i,\bm{\omega}}(T^{\gamma_i+\bm{\omega}\cdot\bm{\delta}}c) \Bigg|^2 \\
        & = \lim_{N \to \infty} \lim_{M \to \infty} \frac{1}{|\Theta_M|} \sum_{\bm{\delta} \in \Theta_M} \Bigg| \frac{1}{|\Theta_N|} \sum_{\bm{\gamma} \in \Theta_N} \prod_{i=1}^k \prod_{\omega_i=1} h_{\bm{\omega}}(T^{\gamma_i-\delta_i+\bm{\omega}\cdot\bm{\delta}}c) \\
        & \qquad \qquad \qquad \qquad \qquad \qquad - \frac{1}{|\Theta_N|} \sum_{\bm{\gamma} \in \Theta_N} \prod_{i=1}^k \prod_{\omega_i=1} h_{\bm{\omega}}(T^{\gamma_i+\bm{\omega}\cdot\bm{\delta}}c) \Bigg|^2.
    \end{align*}
    Unique ergodicity of orbit-closures (\Cref{thm:minimal finite order systems are u.e.}) justifies interchanging the limits, and then we can perform a change of variables $\bm{\gamma}' \to \bm{\gamma} - \bm{\delta}$ in the first term to see that the limit is zero.
\end{proof}

\begin{theorem}
    $\rho_a$ is an Erd\H{o}s cube measure.
\end{theorem}

\begin{proof}
    Let $f_{\bm{\omega}} \colon X \to [0,1]$ be a continuous function for each $\bm{\omega} \in \bracket{k} \setminus \{\bm{0}\}$ such that
    \begin{equation*}
        \int_{X^{\bracket{k}}} 1 \otimes \bigotimes_{\bm{\omega}\ne \bm{0}} f_{\bm{\omega}}~d\rho_a > 0.
    \end{equation*}
    For notational convenience, let $f_{\bm{0}} = 1$ and $f_{\bm{v}} = 1$ for $\bm{v} \in \Z^k \setminus \bracket{k}$.
    Let $\Theta$ be a F{\o}lner sequence in $\Gamma^k$.
    Then by \Cref{lem:equivalent to face action},
    \begin{align*}
        \lim_{N \to \infty} \frac{1}{|\Theta_N|} & \sum_{\bm{\gamma} \in \Theta_N} \int_{X^{\bracket{k}}} \left( \bigotimes_{\bm{\omega} \in \bracket{k}} f_{\bm{\omega}} \right) \cdot \prod_{i=1}^k \bigotimes_{\bm{\omega} \in \bracket{k}} T^{\gamma_i} f_{\bm{\omega}+e_i}~d\rho_a \\
        & = \lim_{N \to \infty} \frac{1}{|\Theta_N|} \sum_{\bm{\gamma} \in \Theta_N} \int_{X^{\bracket{k}}} \left( \bigotimes_{\bm{\omega} \in \bracket{k}} f_{\bm{\omega}} \right) \cdot T_i^{\gamma_i} \prod_{i=1}^k \bigotimes_{\bm{\omega} \in \bracket{k}} f_{i,\bm{\omega}}~d\rho_a \\
        & \geq \lim_{N \to \infty} \frac{1}{|\Theta_N|} \sum_{\bm{\gamma} \in \Theta_N} \int_{X^{\bracket{k}}} \left( \bigotimes_{\bm{\omega} \in \bracket{k}} f_{\bm{\omega}} \right) \cdot T_i^{\gamma_i} \prod_{i=1}^k \bigotimes_{\bm{\omega} \in \bracket{k}} f_{\bm{\omega}}~d\rho_a
    \end{align*}
    Now, $\rho_a$ is invariant under the face actions by construction, so applying the mean ergodic theorem, this final quantity is equal to
    \begin{equation*}
        \int_{X^{\bracket{k}}} F \cdot \E[F \mid \I_1] \cdot \ldots \cdot \E[F \mid \I_k]~d\rho_a,
    \end{equation*}
    where $F = \bigotimes_{\bm{\omega} \in \bracket{k}} f_{\bm{\omega}}$ and $\I_i$ is the $\sigma$-algebra of $T_i$-invariant sets.
    Finally, by Chu's inequality (\cite[Lemma 1.6]{chu}),
    \begin{equation*}
        \int_{X^{\bracket{k}}} F \cdot \E[F \mid \I_1] \cdot \ldots \cdot \E[F \mid \I_k]~d\rho_a \geq \left( \int_{X^{\bracket{k}}} F~d\rho_a \right)^{k+1} > 0.
    \end{equation*}
    This proves that $\rho_a$ is an Erd\H{o}s cube measure.
\end{proof}

\section{Proof of \Cref{cor:main}} \label{sec:proof of corollary}

In this short section, we deduce \Cref{cor:main} from \Cref{thm-main1}.
The small subtlety in our deduction is that $(\N, \times)$ is only a semigroup rather than a group.
However, we can still obtain \Cref{cor:main} from \Cref{thm-main1} with the aid of the following elementary observation:

\begin{lemma} \label{lem:dilate to integers}
    Suppose $k \in \N$ and $B_1, \ldots, B_k \subseteq \Q_{>0}$ are nonempty sets such that
    \begin{equation*}
        B_1 \cdot \ldots \cdot B_k \subseteq \N.
    \end{equation*}
    Then there exist numbers $\lambda_1, \ldots, \lambda_k \in \Q_{>0}$ with $\prod_{j=1}^k \lambda_j = 1$ such that $\lambda_jB_j \subseteq \N$ for every $j \in \{1, \ldots, k\}$.
    Hence, taking $B'_j = \lambda_j B_j$, there are sets $B'_1, \ldots, B'_k \subseteq \N$ with $|B'_j| = |B_j|$ such that
    \begin{equation*}
        B'_1 \cdot \ldots \cdot B'_k = B_1 \cdot \ldots \cdot B_k.
    \end{equation*}
\end{lemma}

\begin{proof}
    We prove the lemma by induction on $k$.
    Let $k = 2$.
    Let $b_1 \in B_1$ be arbitrary.
    Then $b_1 B_2 \subseteq \N$.
    Let $\tilde{B}_1 = B_1/b_1$ and $\tilde{B}_2 = b_1B_2 \subseteq \N$.
    Now let $g = \gcd(\tilde{B}_2)$.
    We can write $g = \sum_{i=1}^r a_i n_i$ for some integers $a_1, \ldots, a_r \in \Z$ and $n_1, \ldots, n_r \in \tilde{B}_2$.
    Therefore,
    \begin{equation*}
        g \tilde{B}_1 = \sum_{i=1}^r a_i \underbrace{(n_i \tilde{B}_1)}_{\subseteq \N} \subseteq \Z.
    \end{equation*}
    Moreover, $g > 0$ and $B_1 \subseteq \Q_{>0}$, so $g \tilde{B}_1 \subseteq \N$.
    We finally take $B'_1 = g \tilde{B}_1$ and $B'_2 = \tilde{B}_2/g$, or equivalently, $\lambda_1 = g/b_1$ and $\lambda_2 = b_1/g$.

    Suppose now that the lemma holds for product sets of size at most $k$ for some $k \ge 2$, and let $B_1, \ldots, B_k, B_{k+1} \subseteq \Q_{>0}$ with $B_1 \cdot \ldots \cdot B_k \cdot B_{k+1} \subseteq \N$.
    Let $C = B_1 \cdot \ldots \cdot B_k$.
    By the $k=2$ case, there exist $\lambda_{k+1} \in \Q_{>0}$ such that $C/\lambda_{k+1}, \lambda_{k+1} B_{k+1} \subseteq \N$.
    Taking $\tilde{B}_1 = B_1/\lambda_{k+1}$, we have $\tilde{B}_1 \cdot B_2 \cdot \ldots \cdot B_k \subseteq \N$, so by the induction hypothesis, there exist numbers $\tilde{\lambda}_1, \lambda_2, \ldots, \lambda_k \in \Q_{>0}$ with $\tilde{\lambda}_1 \cdot \prod_{j=2}^k \lambda_j = 1$ such that $\tilde{\lambda}_1 \tilde{B}_1, \lambda_2 B_2, \ldots, \lambda_k B_k \subseteq \N$.
    Taking $\lambda_1 = \tilde{\lambda}_1/\lambda_{k+1}$, we then have $\tilde{\lambda}_1 \tilde{B}_1 = \lambda_1 B_1$, so $\lambda_j B_j \subseteq \N$ for every $j \in \{1, \ldots, k, k+1\}$, and
    \begin{equation*}
        \prod_{j=1}^{k+1} \lambda_j = \frac{\tilde{\lambda}_1}{\lambda_{k+1}} \cdot \prod_{j=2}^k \lambda_j \cdot \lambda_{k+1} = 1.
    \end{equation*}
\end{proof}

We now prove \Cref{cor:main}:

\begin{proof}[Proof of \Cref{cor:main}]
    Let $A \subseteq \N$ be a set of positive multiplicative upper Banach density, and let $k \in \N$.
    The set $A$ also has positive upper Banach density when viewed as a subset of the abelian group $(\Q_{>0},\times)$ (see, e.g., \cite[Proposition 8.4]{ABS}).
    Hence, by \Cref{thm-main1}, there are infinite sets $B_1, \ldots, B_k \subseteq \Q_{>0}$ such that $B_1 \cdot \ldots \cdot B_k \subseteq A$.
    Applying \Cref{lem:dilate to integers}, there then exist infinite sets $B'_1, \ldots, B'_k \subseteq \N$ such that $B'_1 \cdot \ldots \cdot B'_k = B_1 \cdot \ldots \cdot B_k \subseteq A$.
    This proves \Cref{cor:main}.
\end{proof}

\end{document}